\documentclass[a4paper,11pt]{article}
\usepackage{amsmath,amssymb,amsthm}
\usepackage{a4wide}
\usepackage[dvipsnames]{xcolor}
\usepackage[unicode, final, colorlinks=true]{hyperref}
\usepackage{textgreek}
\usepackage{enumitem, booktabs}
\usepackage{xcolor}
\usepackage{siunitx}
\usepackage{bbm}
\usepackage{mathtools}
\usepackage{tikz}
\usepackage{pgfplots}

\newtoggle{usetikz}
\togglefalse{usetikz}

\pgfplotsset{compat=1.18}
\usetikzlibrary{calc}
\usepgfplotslibrary{groupplots}
\usepgfplotslibrary{colorbrewer}
\usetikzlibrary{external}
\newcommand{\R}{\mathbb{R}}

\newcommand{\defeq}{\coloneqq}

\newtheorem{definition}{Definition}[section]
\newtheorem{theorem}{Theorem}[section]

\newtheorem{proposition}[theorem]{Proposition}
\newtheorem{corollary}[theorem]{Corollary}
\newtheorem{remark}[theorem]{Remark}

\begin{document}
\date{
  \small
  $^1$Institute of Geometry and Applied Mathematics,\\ RWTH Aachen University, Im Süsterfeld 2,\\ 52072 Aachen, Germany\\
  \smallskip
  $^2$Department of Mathematical Sciences, \\
   Ulsan National Institute of Science and Technology,\\
   50 UNIST-gil, Eonyang-eup, Ulju-gun,\\
   Ulsan 44919, Republic of Korea\\
  \bigskip
  October 2026
}

\title{A posteriori error estimates for a finite volume approximation of a local sensing Keller--Segel model}
\author{Niklas Kolbe$^1$ \and Kiwoong Kwon$^{2,}$\footnote{Corresponding author. E-mail: \tt{kwkwon@unist.ac.kr}}}

\maketitle

\begin{abstract}
A local sensing Keller--Segel model, in which cells adjust their motility to the local chemoattractant rather than its gradient, and a semi-implicit finite volume scheme for its numerical approximation are studied. Computable a posteriori bounds for the approximation error measured in terms of the $L^\infty(0, T; H^{-1}(\Omega))$ and $L^2(0, T; L^2(\Omega))$ norms of the cell density are derived based on a new stability estimate and suitable reconstructions of the flux and the flux potential of the numerical solution. Optimal scaling of the estimates is numerically demonstrated in 1D and 2D experiments. Moreover, an adaptive mesh refinement method steered by the estimator is proposed for the 1D case, and its efficiency is numerically verified.
\end{abstract}



\section{Introduction}

Chemotaxis is the directed movement of cells in response to chemical signals.
When cells produce a chemical that attracts other cells, the interaction between
cell movement and chemical production can lead to aggregation and spatial
pattern formation. The Keller--Segel models describe this process through coupled
partial differential equations for the cell density and the chemical
concentration \cite{1970KellerSegel,1971KellerSegel}. In the common gradient sensing form, cell movement combines diffusion with a chemotactic drift driven by the gradient of the chemical
concentration. Widely studied variants include the minimal model, with constant
diffusion and chemotactic sensitivity coefficients, and the classical logarithmic model,
with constant diffusion and a sensitivity inversely proportional to the chemical
concentration \cite[equations~(7)--(8)]{2017YoonKim}. Depending on the model, spatial
dimension, and initial data, chemotactic aggregation can lead to finite-time
blow-up \cite{1997HerreroVelazquez}. This possibility motivates the study of models
that admit global solutions while retaining the ability to describe cell
aggregation and spatial pattern formation~\cite{2017YoonKim,2019Desvillettes}.

Local sensing models retain a link between diffusion and chemotactic sensitivity
inherited from the original random-walk derivation of the Keller--Segel equations
\cite{1971KellerSegel,2017YoonKim}. In that derivation, the cell motility (a function of the chemical concentration) connects the two processes.
The local sensing formulation describes cells that respond to the local chemical
concentration by changing their motility, without requiring them to measure a
chemical gradient. The cell flux then takes the form $J=-\nabla(\gamma(v)\rho)$,
where $\gamma(v)$ denotes the motility, and its diffusion and chemotactic drift
components are determined by the same function. A decreasing motility promotes
accumulation in regions of high chemical concentration. This mechanism was
studied in nutrient-consumption models exhibiting traveling waves \cite{2015YoonKim}
and in models of cell aggregation driven by a self-produced chemical
\cite{2017YoonKim,2019Desvillettes}. The latter provide a framework for studying
global existence and aggregation within the same model. Related experiments and
kinetic modeling also demonstrate stripe formation through density-suppressed
motility \cite{2011LiuEtAl,2012FuEtAl}, and a rigorous mean-field derivation connects a
signal-dependent model to an interacting particle system \cite{2026BolChenLi}.

For the parabolic--elliptic system, global existence results cover general classes
of positive motilities, including functions that decay to zero
\cite{2021LiJiang,2021JiangLaurencot}. Global solutions can nevertheless develop
increasingly concentrated density profiles. For exponential motility, the
two-dimensional parabolic--parabolic model exhibits critical-mass behavior, with
uniformly bounded solutions for subcritical mass and unbounded growth as time
tends to infinity for suitable supercritical initial data
\cite{2020JinWang,2021BurgerLaurencotTrescases}.

The Laplace form of the density equation $\partial_t\rho=\Delta(\gamma(v)\rho)$
admits a duality estimate: testing the equation with the inverse Neumann
Laplacian of the deviation of the density from its mean controls the density in
the dual space of $H^1$ and the integral of $\gamma(v)\rho^2$ over space and
time \cite{2021BurgerLaurencotTrescases}. Global weak solutions were constructed
from this estimate and the entropy structure, and the density equation was
interpreted in a very weak sense, with both spatial derivatives of the diffusion
term on the test function, as for the porous medium equation
\cite{2021BurgerLaurencotTrescases}; see also \cite{2019Desvillettes}. This
formulation requires no spatial derivative of the density.

For gradient sensing Keller--Segel models, convergence was established for a
finite volume scheme when the initial mass lies below a threshold
\cite{2006Filbet}, and error estimates were obtained for a conservative upwind
finite element scheme \cite{2007Saito} and for interior penalty discontinuous
Galerkin methods \cite{2008EpshteynKurganov} under regularity assumptions on
the exact solution. Numerical analysis of local sensing chemotaxis models is
more limited. For a parabolic--elliptic model with logistic growth, generalized
finite difference schemes were shown to converge under a restriction on the
time step \cite{2023HerreroHervasNegreanuVargas,2024HerreroHervas}. For the
local sensing model without logistic growth and with exponential motility, a
linearly implicit finite volume scheme was proposed in
\cite{2025HerdaTrescasesZurek}. This scheme preserves nonnegativity and mass,
and discrete entropy and duality estimates yield unconditional stability and
subsequential convergence to weak solutions.

Numerical experiments in \cite{2025HerdaTrescasesZurek} reveal difficulties in
resolving the long-time concentration of the density. Nonnegativity and mass
conservation impose a mesh-dependent upper bound on the discrete density, whose
maximum may saturate on coarse meshes. The slow growth of this maximum also
makes infinite-time blow-up difficult to distinguish from bounded growth. These
observations motivate further numerical investigation using refined meshes and
adaptive strategies. A computable a posteriori error estimate can quantify the
error of a given numerical solution and guide mesh refinement in regions that
require greater resolution. Such a posteriori error estimates were recently derived for the classical parabolic--elliptic Keller--Segel
model, using a discontinuous Galerkin
approximation \cite{2024GiesselmannKwon} and a combined finite volume and
finite element approximation \cite{2025HoffmannGiesselmann}
and for a finite volume approximation of a gradient sensing
Keller--Segel model with power-law diffusion \cite{2025GiesselmannKolbe}.

In this work, we derive residual-based a posteriori error estimates for a
fully discrete finite volume approximation of the parabolic--elliptic
local sensing Keller--Segel system, namely the quasi-stationary adaptation of
the scheme \cite{2025HerdaTrescasesZurek}.
We expect our approach to extend to the parabolic--parabolic system studied in
\cite{2025HerdaTrescasesZurek}. Since this extension involves additional
technical difficulties, we restrict ourselves to the parabolic--elliptic system
and leave the parabolic--parabolic case to future work. We consider positive, bounded, and
locally Lipschitz continuous motilities, which include motilities that decay to
zero such as the exponential motility. 
The nonlinear flux potential $m=\gamma(v)\rho$ couples the errors of the density and
the chemical signal, giving rise to a mixed term in the error analysis. 
 As the diffusion coefficient $\gamma(v)$ is not bounded away from zero, we first derive an $L^\infty$ bound for the chemoattractant
concentration of the weak solution on a finite time interval $(0,T)$ by a
comparison argument \cite{2020FujieJiang,2021LiJiang}. This bound
confines the values of the concentration to a bounded
interval, on which the motility is Lipschitz continuous and bounded from below.
The very weak formulation for the local sensing model does not allow for stability estimates within an $L^2$ framework, such as the ones in \cite{2024GiesselmannKwon} and \cite[Section~2.1]{2025GiesselmannKolbe}. 
We therefore measure the density error in $H^{-1}$, in analogy
with the duality estimate and with the stability estimate for power-law
diffusion in \cite[Section~2.2]{2025GiesselmannKolbe}. We then establish a
stability estimate for sufficiently regular
approximations that controls the density error in terms of the initial error
and residuals. The a posteriori analyses of the gradient sensing
models
\cite{2024GiesselmannKwon,2025GiesselmannKolbe,2025HoffmannGiesselmann} rely on
conditional stability estimates. In contrast, our stability estimate for the
local sensing system requires no a posteriori smallness condition.
Our stability estimate is applied to reconstructions of the density and the flux that
satisfy the continuity equation exactly, an $H^1$-conforming flux potential
reconstruction, and an elliptic reconstruction that satisfies the elliptic
equation exactly.
The chosen reconstructions remain stable under mesh refinement on finite volumes and therefore enable the design of an adaptive method.
By bounding the resulting residuals in terms of the numerical
solution and the problem data, we obtain a computable upper bound for the
density error.

Numerical experiments in one and two space dimensions show that the decay of the estimator components is consistent with the scaling of the error in a piecewise-constant first order approximation. 
We also split the estimator into local contributions to enable an adaptive
mesh refinement method in one space dimension. We demonstrate that the adaptive method outperforms computations on uniform meshes in terms of the estimator and also improves efficiency if sufficiently large refinement levels are allowed.

We organize the remainder of the paper as follows. In Section 2, we introduce the
model, define weak solutions, and derive an $L^\infty$ bound for the
chemoattractant concentration. In Section 3, we establish the stability
estimate. In Section 4, we present the finite volume scheme for the
parabolic--elliptic model. In Section 5, we construct the density, flux potential,
flux, and elliptic reconstructions. In Section 6, we estimate the elliptic
reconstruction error and combine the residual bounds with the stability estimate
to obtain a computable a posteriori error bound. In Section 7, we examine the
convergence of the estimator components and test an adaptive mesh refinement
method based on their local contributions before we conclude the article in Section 8.

\section{Model}
In this section, we introduce the local sensing Keller--Segel model, define
weak solutions, and obtain an $L^\infty$ bound for the chemoattractant
concentration.

\subsection{Local sensing Keller--Segel model}\label{subsec:model}

Let $\Omega\subset\R^d$ for dimension $d\in\{1,2, 3\}$ be a bounded convex polyhedral domain with outward unit normal $\nu$. We consider the system
\begin{equation}\label{eq:model}
\left\{  \begin{aligned}
    \partial_t\rho=\Delta\bigl(\gamma(v)\rho\bigr), \\
    -\Delta v+\beta v=\alpha\rho,
\end{aligned} \right. \qquad (t,x) \in (0,T) \times \Omega,
\end{equation}
where $\rho$ denotes the cell density, $v$ denotes the chemoattractant concentration, and $\gamma$ is the motility function. The positive constants $\alpha$ and $\beta$ denote the production and decay rates of the chemoattractant, respectively. System \eqref{eq:model} is subject to the homogeneous Neumann boundary conditions
\[
\nabla(\gamma(v)\rho)\cdot\nu=0,
\qquad
\nabla v\cdot\nu=0
\qquad\text{on }(0,T)\times\partial\Omega,
\]
under which neither cells nor chemoattractant cross the boundary, and to the initial condition
\[
\rho(0,\cdot)=\rho^0
\qquad\text{in }\Omega,
\]
where the initial density $\rho^0\in L^2(\Omega)$ is nonnegative. The motility is a scalar function of the attractant concentration. We assume that
\begin{equation}\label{eq:gamma-assumptions}
\gamma \in C^{0,1}_{\mathrm{loc}}([0,\infty)), \qquad 0 < \gamma(s) \le \gamma_1 \quad \text{for all } s \ge 0.
\end{equation}
On each interval $[0,B]$ with $B>0$, the motility is Lipschitz continuous and bounded from below by a positive constant.

\begin{remark}[Examples of motility functions]
Assumption \eqref{eq:gamma-assumptions} is satisfied by the exponential motility $\gamma(s)=e^{-\chi s}$ with $\chi>0$, which we use in the numerical experiments in Section~\ref{sec:experiments}, and by the algebraically decaying motility $\gamma(s)=(1+s)^{-k}$ with $k>0$. In both cases, $\gamma_1=1$. Power-law motilities $\gamma(s)=s^{-k}$ with $k>0$ are singular at $s=0$ and do not satisfy \eqref{eq:gamma-assumptions}.
\end{remark}

Next, we define weak solutions to the local sensing Keller--Segel model.
\begin{definition}[Weak solution to \eqref{eq:model}]\label{def:weak-solution}
Let $T>0$ and let $\rho^0 \in L^2(\Omega)$ be nonnegative. We call a pair $(\rho,v)$ a weak solution to \eqref{eq:model} on $(0,T)$ if
\[
\rho \ge 0, \qquad v \ge 0 \qquad \text{a.e. on } (0,T)\times\Omega,
\]
\[
\rho \in L^2((0,T)\times\Omega), \qquad
\partial_t\rho \in L^2(0,T;(H^1(\Omega))'), \qquad
v \in L^\infty(0,T;H^1(\Omega)),
\]
and if the following equations hold:
\begin{equation}\label{eq:weak-rho}
\int_0^T\!\!\int_{\Omega} \rho\,\partial_t\phi \, dx\,dt
+ \int_{\Omega}\rho^0(x)\phi(x,0)\,dx
+ \int_0^T\!\!\int_{\Omega}\gamma(v)\rho\,\Delta\phi \, dx\,dt = 0
\end{equation}
for every $\phi \in C^1([0,T];H^2(\Omega))$ satisfying $\phi(\cdot,T)=0$ and $\partial_\nu\phi=0$ on $\partial\Omega\times[0,T]$, and
\begin{equation}\label{eq:weak-v}
\int_0^T\!\!\int_{\Omega} \nabla v \cdot \nabla \psi \, dx\,dt
+ \beta \int_0^T\!\!\int_{\Omega} v\psi \, dx\,dt
- \alpha \int_0^T\!\!\int_{\Omega} \rho\psi \, dx\,dt = 0
\end{equation}
for every $\psi \in C^\infty(\overline\Omega \times [0,T])$.
\end{definition}

\begin{remark}[Regularity of weak solutions]
In Definition \ref{def:weak-solution}, we adapt the distributional formulation
in \cite{2025HerdaTrescasesZurek} to the present parabolic--elliptic system.
Compared with that formulation, Definition \ref{def:weak-solution} additionally
requires the temporal regularity
\[
\partial_t\rho\in L^2(0,T;(H^1(\Omega))'),
\]
which we use in the stability estimate in Theorem \ref{thm:hminus1-stability}.
\end{remark}

\subsection{Bound for the chemoattractant concentration}\label{subsec:signal-bound}
For motilities that decay to zero, such as the exponential motility, the diffusion coefficient $\gamma(v)$ in \eqref{eq:model} is bounded away from zero only if $v$ is bounded. In this subsection, we obtain an $L^\infty$ bound, depending on $T$, for the chemoattractant concentration of every weak solution on $(0,T)$ and the resulting positive lower bound for the motility. Both bounds enter the stability estimate in Theorem~\ref{thm:hminus1-stability}. The bound follows from a comparison argument developed in \cite{2020FujieJiang,2021LiJiang}; see also \cite[Lemma~3.1]{2021JiangLaurencot}.

For $f\in L^2(\Omega)$, let $A_\beta^{-1}f\in H^1(\Omega)$ denote the unique solution $u$ of
\begin{equation}\label{eq:resolvent-weak}
\int_\Omega \nabla u\cdot\nabla\eta\,dx
+\beta\int_\Omega u\eta\,dx
=\int_\Omega f\eta\,dx
\qquad\text{for every }\eta\in H^1(\Omega).
\end{equation}
The operator $A_\beta^{-1}$ is bounded and symmetric on $L^2(\Omega)$. Choosing $\eta=1$ in \eqref{eq:resolvent-weak} shows that $f-\beta u$ has zero mean. Since $\Omega$ is convex, elliptic regularity for the homogeneous Neumann problem \cite[Theorem~2.4.3 and Remark~2.6.9]{1992Grisvard} gives
\begin{equation}\label{eq:resolvent-regularity}
u\in H^2(\Omega),
\qquad
-\Delta u+\beta u=f\quad\text{in }\Omega,
\qquad
\partial_\nu u=0\quad\text{on }\partial\Omega.
\end{equation}
If $f\ge0$ almost everywhere in $\Omega$, then choosing $\eta=\min\{u,0\}$ in \eqref{eq:resolvent-weak} gives $u\ge0$ almost everywhere in $\Omega$. Consequently, $A_\beta^{-1}$ preserves the order of functions in $L^2(\Omega)$.

\begin{proposition}[Bound for the chemoattractant concentration]\label{prop:signal-bound}
Let $(\rho,v)$ be a weak solution to \eqref{eq:model} in the sense of Definition \ref{def:weak-solution}, and set
\[
v^0\defeq\alpha A_\beta^{-1}\rho^0,
\qquad
B_0\defeq\|v^0\|_{L^\infty(\Omega)}.
\]
Then $B_0<\infty$, $v\in H^1(0,T;L^2(\Omega))$ with $v(0)=v^0$ and the chemoattractant satisfies
\begin{equation}\label{eq:resolvent-identity}
\partial_t v+\alpha\gamma(v)\rho
=\alpha\beta A_\beta^{-1}\bigl(\gamma(v)\rho\bigr)
\qquad\text{a.e. in }(0,T)\times\Omega.
\end{equation}
Moreover, for every $t\in[0,T]$ it holds
\begin{equation}\label{eq:signal-bound}
0\le v(t,x)\le e^{\beta\gamma_1 t}\,v^0(x)
\qquad\text{for a.e. }x\in\Omega,
\end{equation}
and, in particular,
\[
0\le v\le B_T\defeq B_0e^{\beta\gamma_1T},
\qquad
\gamma(v)\ge\gamma_{0,T}\defeq\min_{0\le s\le B_T}\gamma(s)>0
\]
almost everywhere in $(0,T)\times\Omega$.
\end{proposition}

\begin{proof}
Set $m\defeq\gamma(v)\rho$. Assumption \eqref{eq:gamma-assumptions} and $\rho\ge0$ imply
\begin{equation}\label{eq:m-bounds}
0\le m\le\gamma_1\rho
\qquad\text{a.e. in }(0,T)\times\Omega.
\end{equation}
By \eqref{eq:resolvent-regularity}, $v^0\in H^2(\Omega)$. Since $d\le3$, the Sobolev embedding $H^2(\Omega)\hookrightarrow L^\infty(\Omega)$ \cite[Theorem~4.12]{2003AdamsFournier} gives $B_0<\infty$.

Equation \eqref{eq:weak-v} gives $v(t)=\alpha A_\beta^{-1}\rho(t)$ for almost every $t\in(0,T)$. Let $\psi\in L^2(\Omega)$ and set $\varphi\defeq\alpha A_\beta^{-1}\psi$. By \eqref{eq:resolvent-regularity}, $\Delta\varphi=\beta\varphi-\alpha\psi$, and the symmetry of $A_\beta^{-1}$ gives, for almost every $t\in(0,T)$,
\[
\int_\Omega\rho(t)\varphi\,dx=\int_\Omega v(t)\psi\,dx,
\qquad
\int_\Omega m(t)\Delta\varphi\,dx
=\int_\Omega\bigl(\alpha\beta A_\beta^{-1}m(t)-\alpha m(t)\bigr)\psi\,dx,
\]
as well as
\[
\int_\Omega\rho^0\varphi\,dx=\int_\Omega v^0\psi\,dx.
\]
Let $\eta\in C^1([0,T])$ satisfy $\eta(T)=0$. By \eqref{eq:resolvent-regularity}, $\phi(x,t)=\eta(t)\varphi(x)$ is an admissible test function in \eqref{eq:weak-rho}, which then reads
\begin{equation}\label{eq:signal-weak-time}
\int_0^T\eta'(t)\int_\Omega v\psi\,dx\,dt
+\eta(0)\int_\Omega v^0\psi\,dx
+\int_0^T\eta(t)\int_\Omega\bigl(\alpha\beta A_\beta^{-1}m-\alpha m\bigr)\psi\,dx\,dt
=0.
\end{equation}
Since $m\in L^2((0,T)\times\Omega)$ and $A_\beta^{-1}$ is bounded on $L^2(\Omega)$, we have $\alpha\beta A_\beta^{-1}m-\alpha m\in L^2((0,T)\times\Omega)$. Identity \eqref{eq:signal-weak-time} holds for every $\psi\in L^2(\Omega)$ and every such $\eta$. Hence $v\in H^1(0,T;L^2(\Omega))$, identity \eqref{eq:resolvent-identity} holds, and $v(0)=v^0$.

For almost every $t\in(0,T)$, estimate \eqref{eq:m-bounds} and the order preservation of $A_\beta^{-1}$ yield
\[
0\le A_\beta^{-1}m(t)
\le\gamma_1A_\beta^{-1}\rho(t)
=\frac{\gamma_1}{\alpha}v(t)
\qquad\text{a.e. in }\Omega.
\]
Using $m\ge0$ in \eqref{eq:resolvent-identity}, we obtain $\partial_t v\le\beta\gamma_1v$ almost everywhere in $(0,T)\times\Omega$. The function $w(t)\defeq e^{-\beta\gamma_1t}v(t)$ therefore satisfies $w\in H^1(0,T;L^2(\Omega))$ and $\partial_t w\le0$ almost everywhere in $(0,T)\times\Omega$. For every $t\in[0,T]$,
\[
w(t)=v^0+\int_0^t\partial_t w(s)\,ds\le v^0
\qquad\text{a.e. in }\Omega,
\]
which is the upper bound in \eqref{eq:signal-bound}. The lower bound follows from $v\ge0$ and the continuity of $v$ in time with values in $L^2(\Omega)$. Finally, \eqref{eq:signal-bound} gives $0\le v\le B_T$ almost everywhere in $(0,T)\times\Omega$, and $\gamma_{0,T}>0$ because $\gamma$ is continuous and positive on $[0,B_T]$.
\end{proof}

\begin{remark}
If $\rho^0\in L^\infty(\Omega)$, due to the order preservation of $A_\beta^{-1}$ and $A_\beta^{-1}1=\beta^{-1}$ it holds
\[
B_0\le\frac{\alpha}{\beta}\|\rho^0\|_{L^\infty(\Omega)}.
\]
\end{remark}

\section{Stability framework}
In this section, a stability framework is established that will be essential for the derivation of a posteriori error estimates for the local sensing Keller--Segel model. In the process, the difference between the PDE solution of~\eqref{eq:model} and a sufficiently regular approximation is analyzed.

The approximation we consider is denoted by
\[
(\bar\rho,\bar v,\bar m,\bar J)
\]
and defined on $(0,T)\times\Omega$ with the regularity assumptions
\begin{equation}\label{eq:strongapprox-regularity}
\left\{  \begin{aligned}
\bar\rho &\in L^\infty((0,T)\times\Omega), \qquad
\partial_t\bar\rho \in L^2(0,T;(H^1(\Omega))'),\\
\bar v,\bar m &\in L^2(0,T;H^1(\Omega)), \qquad
                \bar J \in L^2(0,T;L^2(\Omega;\R^d)).
  \end{aligned} \right.
\end{equation}
We assume that $\bar v\ge0$ almost everywhere on $(0,T)\times\Omega$, so that $\gamma(\bar v)$ is well-defined. For almost every $t\in(0,T)$, the approximations $\bar \rho$ and $\bar v$ satisfy the continuity equation and the elliptic equation in \eqref{eq:model} in the sense that
\begin{equation}\label{eq:recon-mass}
\left\langle\partial_t\bar\rho,\xi\right\rangle_{(H^1(\Omega))',H^1(\Omega)}
-\int_\Omega\bar J\cdot\nabla\xi\,dx=0
\qquad\text{for every }\xi\in H^1(\Omega),
\end{equation}
and
\begin{equation}\label{eq:recon-elliptic}
\int_\Omega\nabla\bar v\cdot\nabla\eta\,dx
+\beta\int_\Omega\bar v\eta\,dx
=\alpha\int_\Omega\bar\rho\eta\,dx
\qquad\text{for every }\eta\in H^1(\Omega).
\end{equation}
The approximations $\bar J$ and $\bar m$ connect equations \eqref{eq:recon-mass} and \eqref{eq:recon-elliptic}. We define the residuals
\begin{align}
R_f &\defeq \bar J + \nabla \bar m,
\label{eq:Rflux} \\
R_c &\defeq \bar m - \gamma(\bar v)\bar\rho,
\label{eq:Rcons}
\end{align}
which measure the deviation from the flux law $J=-\nabla m$ and from the constitutive law $m=\gamma(v)\rho$, respectively. The regularity assumptions \eqref{eq:strongapprox-regularity} and the boundedness of $\gamma$ give
\begin{equation}\label{eq:residuals}
R_f \in L^2(0,T;L^2(\Omega;\R^d)),
\qquad
R_c \in L^2(0,T;L^2(\Omega)).
\end{equation}

In the following result, we use the zero-average spaces
\[
L^2_\diamond(\Omega) \defeq \left\{f \in L^2(\Omega) \,\middle|\, \int_{\Omega}f\,dx = 0 \right\}, \quad
H^1_\diamond(\Omega) \defeq \left\{\phi \in H^1(\Omega) \,\middle|\, \int_{\Omega}\phi\,dx = 0 \right\},
\]
and define
\[
H^{-1}_\diamond(\Omega) \defeq
\left\{F \in (H^1(\Omega))' \,\middle|\,
\langle F,1\rangle_{(H^1(\Omega))',H^1(\Omega)} = 0 \right\}.
\]
For $F\in H^{-1}_\diamond(\Omega)$, we set
\[
\|F\|_{H^{-1}(\Omega)} \defeq
\sup_{\phi \in H^1_\diamond(\Omega)\setminus\{0\}}
\frac{\langle F,\phi\rangle_{(H^1(\Omega))',H^1(\Omega)}}
{\|\nabla \phi\|_{L^2(\Omega)}}.
\]
We next state and prove the central stability result of this work.

\begin{theorem}\label{thm:hminus1-stability}
Let $(\rho,v)$ be a weak solution of \eqref{eq:model} in the sense of Definition \ref{def:weak-solution}, and let $(\bar\rho,\bar v,\bar m,\bar J)$ be an approximation on $(0,T)\times\Omega$ satisfying \eqref{eq:strongapprox-regularity}--\eqref{eq:residuals}.
Assume further that
\[
\rho^0-\bar\rho(0)\in H^{-1}_\diamond(\Omega)
\]
and define the error
\[
e \defeq \rho-\bar\rho.
\]
Let $B_T$ and $\gamma_{0,T}$ be the constants in Proposition \ref{prop:signal-bound}, set
\begin{equation}\label{eq:explicit-constants}
M_T \defeq \|\bar\rho\|_{L^\infty((0,T)\times\Omega)},
\qquad
\tilde \alpha \defeq \frac{\alpha}{\min\{\beta, 1\}},
\end{equation}
and let $L_{\gamma,T}$ denote the Lipschitz constant of $\gamma$ on the interval $[0,\max\{B_T,\alpha M_T/\beta\}]$. Define
\begin{equation}\label{eq:lambda-constant}
\Lambda_T \defeq 1 + \frac{2L_{\gamma,T}^2 M_T^2 \tilde \alpha^2}{\gamma_{0,T}}.
\end{equation}
Then the error satisfies
\begin{multline}\label{eq:gronwall-estimate}
\sup_{0 \le t \le T}\|e(t)\|_{H^{-1}(\Omega)}^2
+ \int_0^T \|e(t)\|_{L^2(\Omega)}^2\,dt \\
\le \left(1+\frac{1}{\gamma_{0,T}}\right)e^{\Lambda_T T}\Biggl(
\|e(0)\|_{H^{-1}(\Omega)}^2
+ \int_0^T \|R_f(t)\|_{L^2(\Omega)}^2\,dt
+ \frac{2}{\gamma_{0,T}}\int_0^T \|R_c(t)\|_{L^2(\Omega)}^2\,dt
\Biggr).
\end{multline}
\end{theorem}

\begin{proof}
Set
\[
m \defeq \gamma(v)\rho,
\qquad
z \defeq v-\bar v,
\qquad
\mu \defeq m-\bar m.
\]
Taking test functions of the form $\phi(x,t)=\eta(t)$ in \eqref{eq:weak-rho} yields conservation of the total mass of $\rho$. Choosing $\xi=1$ in \eqref{eq:recon-mass} yields conservation of the total mass of $\bar\rho$. Since $\rho^0-\bar\rho(0)\in H^{-1}_\diamond(\Omega)$,
\[
\int_\Omega\bigl(\rho^0-\bar\rho(0)\bigr)\,dx=0.
\]
Combining this identity with conservation of both masses gives
\[
\int_{\Omega} e(t,x)\,dx = 0
\qquad \text{for a.e. } t \in (0,T).
\]
Subtracting \eqref{eq:recon-elliptic} from \eqref{eq:weak-v}, we obtain, for almost every $t\in(0,T)$,
\begin{equation}\label{eq:error-elliptic}
\int_\Omega \nabla z\cdot\nabla\eta\,dx
+\beta\int_\Omega z\eta\,dx
=\alpha\int_\Omega e\eta\,dx
\qquad\text{for every }\eta\in H^1(\Omega).
\end{equation}
For every $\psi\in H^2(\Omega)$ satisfying $\partial_\nu\psi=0$ on $\partial\Omega$, equations \eqref{eq:weak-rho}, \eqref{eq:recon-mass}, and \eqref{eq:Rflux} give, for almost every $t\in(0,T)$,
\begin{equation}\label{eq:error-mass-rewritten}
\left\langle\partial_t e,\psi\right\rangle_{(H^1(\Omega))',H^1(\Omega)}
=
\int_\Omega\mu\,\Delta\psi\,dx
-\int_\Omega R_f\cdot\nabla\psi\,dx.
\end{equation}
For almost every $t\in(0,T)$, let $\phi(t)\in H^1_\diamond(\Omega)$ be the unique weak solution of
\[
\int_\Omega\nabla\phi(t)\cdot\nabla\eta\,dx
=
\int_\Omega e(t)\eta\,dx
\qquad
\text{for every }\eta\in H^1_\diamond(\Omega).
\]
Since $\Omega$ is convex and $e(t)\in L^2_\diamond(\Omega)$, elliptic regularity for the homogeneous Neumann problem \cite[Theorem~2.4.3 and Remark~2.6.9]{1992Grisvard} gives
\[
\phi(t)\in H^2(\Omega),\qquad
-\Delta\phi(t)=e(t)\quad\text{in }\Omega,\qquad
\partial_\nu\phi(t)=0\quad\text{on }\partial\Omega.
\]
The Hilbert-space chain rule and the definition of the $H^{-1}(\Omega)$-norm give
\[
\frac12\frac{d}{dt}\|e(t)\|_{H^{-1}(\Omega)}^2
=
\left\langle\partial_t e(t),\phi(t)\right\rangle_{(H^1(\Omega))',H^1(\Omega)},
\qquad
\|\nabla\phi(t)\|_{L^2(\Omega)}
=
\|e(t)\|_{H^{-1}(\Omega)}.
\]
Taking $\psi=\phi(t)$ in \eqref{eq:error-mass-rewritten} and using $\Delta\phi(t)=-e(t)$, we obtain
\begin{equation}\label{eq:basic-identity}
\frac12\frac{d}{dt}\|e(t)\|_{H^{-1}(\Omega)}^2
+\int_\Omega\mu(t)e(t)\,dx
=
-\int_\Omega R_f(t)\cdot\nabla\phi(t)\,dx
\qquad
\text{for a.e. }t\in(0,T).
\end{equation}

Using the definitions of $e$ and $\mu$ together with \eqref{eq:Rcons}, we obtain
\[
\mu = \gamma(v)e + (\gamma(v)-\gamma(\bar v))\bar\rho - R_c.
\]
Multiplying by $e$ and integrating over $\Omega$ gives
\begin{equation}\label{eq:mu-e-splitting}
\int_\Omega \mu e\,dx
= \int_\Omega \gamma(v)e^2\,dx
+ \int_\Omega (\gamma(v)-\gamma(\bar v))\bar\rho\,e\,dx
- \int_\Omega R_c e\,dx.
\end{equation}
Proposition \ref{prop:signal-bound} gives $\gamma(v)\ge\gamma_{0,T}$ almost everywhere in $(0,T)\times\Omega$, and therefore
\begin{equation}\label{eq:coercive}
\int_\Omega \gamma(v)e^2\,dx \ge \gamma_{0,T} \|e\|_{L^2(\Omega)}^2.
\end{equation}
By Proposition \ref{prop:signal-bound}, $0\le v\le B_T$ almost everywhere in $(0,T)\times\Omega$. Since $|\bar\rho|\le M_T$, the weak maximum principle for \eqref{eq:recon-elliptic} and $\bar v\ge0$ give $0\le\bar v\le\alpha M_T/\beta$. The Lipschitz continuity of $\gamma$ on $[0,\max\{B_T,\alpha M_T/\beta\}]$, the definition of $M_T$, and H\"older's inequality therefore give
\[
\left|\int_\Omega (\gamma(v)-\gamma(\bar v))\bar\rho\,e\,dx\right|
\le L_{\gamma,T} M_T \|z\|_{L^2(\Omega)} \|e\|_{L^2(\Omega)}.
\]
Taking $\eta=z$ in \eqref{eq:error-elliptic} gives
\[
\|\nabla z\|_{L^2(\Omega)}^2
+ \beta \|z\|_{L^2(\Omega)}^2
= \alpha \int_\Omega ez\,dx.
\]
Moreover,
\[
\min\{\beta,1\}\|z\|_{H^1(\Omega)}^2
\le
\|\nabla z\|_{L^2(\Omega)}^2
+ \beta \|z\|_{L^2(\Omega)}^2.
\]
Because $e$ has zero average, the definition of the $H^{-1}(\Omega)$-norm gives
\[
\left|\int_\Omega ez\,dx\right|
\le \|e\|_{H^{-1}(\Omega)}\|\nabla z\|_{L^2(\Omega)}
\le \|e\|_{H^{-1}(\Omega)}\|z\|_{H^1(\Omega)}.
\]
Using $\tilde\alpha$ defined in \eqref{eq:explicit-constants}, we obtain
\[
\|z\|_{L^2(\Omega)}
\le \|z\|_{H^1(\Omega)}
\le \tilde\alpha \|e\|_{H^{-1}(\Omega)}.
\]
Consequently,
\[
\left|\int_\Omega (\gamma(v)-\gamma(\bar v))\bar\rho\,e\,dx\right|
\le L_{\gamma,T} M_T\tilde\alpha
\|e\|_{H^{-1}(\Omega)} \|e\|_{L^2(\Omega)}.
\]
Young's inequality then gives
\begin{equation}\label{eq:mixed-bound}
\left|\int_\Omega (\gamma(v)-\gamma(\bar v))\bar\rho\,e\,dx\right|
\le \frac{\gamma_{0,T}}{4}\|e\|_{L^2(\Omega)}^2
+ \frac{L_{\gamma,T}^2 M_T^2 \tilde\alpha^2}{\gamma_{0,T}}
\|e\|_{H^{-1}(\Omega)}^2.
\end{equation}

The Cauchy--Schwarz and Young inequalities give
\begin{equation}\label{eq:Rcons-bound}
\left|\int_\Omega R_c e\,dx\right|
\le \|R_c\|_{L^2(\Omega)} \|e\|_{L^2(\Omega)}
\le \frac{\gamma_{0,T}}{4}\|e\|_{L^2(\Omega)}^2
+ \frac{1}{\gamma_{0,T}}\|R_c\|_{L^2(\Omega)}^2.
\end{equation}
Combining identity \eqref{eq:mu-e-splitting} with estimates \eqref{eq:coercive}, \eqref{eq:mixed-bound}, and \eqref{eq:Rcons-bound}, we obtain
\begin{equation}\label{eq:mu-lower-bound}
\int_\Omega \mu e\,dx
\ge \frac{\gamma_{0,T}}{2}\|e\|_{L^2(\Omega)}^2
- \frac{L_{\gamma,T}^2 M_T^2 \tilde\alpha^2}{\gamma_{0,T}}
\|e\|_{H^{-1}(\Omega)}^2
- \frac{1}{\gamma_{0,T}}\|R_c\|_{L^2(\Omega)}^2.
\end{equation}

The identity
$\|\nabla\phi\|_{L^2(\Omega)}=\|e\|_{H^{-1}(\Omega)}$
and the Cauchy--Schwarz and Young inequalities give
\begin{equation}\label{eq:Rflux-bound}
\begin{aligned}
\left|\int_\Omega R_f\cdot\nabla\phi\,dx\right|
&\le \|R_f\|_{L^2(\Omega)}\|\nabla\phi\|_{L^2(\Omega)}\\
&=\|R_f\|_{L^2(\Omega)}\|e\|_{H^{-1}(\Omega)}\\
&\le \frac12\|R_f\|_{L^2(\Omega)}^2
+\frac12\|e\|_{H^{-1}(\Omega)}^2.
\end{aligned}
\end{equation}

Substituting estimates \eqref{eq:mu-lower-bound} and
\eqref{eq:Rflux-bound} into identity \eqref{eq:basic-identity},
multiplying by $2$, and using the definition of $\Lambda_T$ in
\eqref{eq:lambda-constant}, we obtain
\begin{equation}\label{eq:hminus1-differential}
\frac{d}{dt}\|e\|_{H^{-1}(\Omega)}^2
+ \gamma_{0,T}\|e\|_{L^2(\Omega)}^2
\le
\Lambda_T\|e\|_{H^{-1}(\Omega)}^2
+ \|R_f\|_{L^2(\Omega)}^2
+ \frac{2}{\gamma_{0,T}}\|R_c\|_{L^2(\Omega)}^2
\end{equation}
for almost every $t \in (0,T)$. Set
\[
y(t) \defeq \|e(t)\|_{H^{-1}(\Omega)}^2,
\qquad
g(t) \defeq \|e(t)\|_{L^2(\Omega)}^2,
\qquad
Y(t) \defeq y(t) + \gamma_{0,T}\int_0^t g(s)\,ds.
\]
Integrating \eqref{eq:hminus1-differential} over $(0,t)$
and using $y(s)\le Y(s)$, we obtain
\[
Y(t)
\le
y(0) + \Lambda_T \int_0^t Y(s)\,ds
+ \int_0^t \|R_f(s)\|_{L^2(\Omega)}^2\,ds
+ \frac{2}{\gamma_{0,T}}\int_0^t \|R_c(s)\|_{L^2(\Omega)}^2\,ds.
\]
Applying Gronwall's inequality gives
\[
Y(t)
\le
e^{\Lambda_T t}\left(
y(0)
+
\int_0^t \|R_f(s)\|_{L^2(\Omega)}^2\,ds
+
\frac{2}{\gamma_{0,T}}\int_0^t \|R_c(s)\|_{L^2(\Omega)}^2\,ds
\right)
\]
for every $t \in [0,T]$. Since $y(t)\le Y(t)$ and
$\gamma_{0,T}\int_0^T g(s)\,ds \le Y(T)$, taking the supremum over
$t\in[0,T]$ proves \eqref{eq:gronwall-estimate}.
\end{proof}

\begin{remark}
Unlike the stability estimates in \cite{2024GiesselmannKwon,2025GiesselmannKolbe}, the estimate in Theorem~\ref{thm:hminus1-stability} requires no a posteriori smallness condition on the initial error or the residuals under its stated assumptions. The existence and stated regularity of the exact solution are assumed, and the constants depend on $T$.
\end{remark}

\begin{remark}[Dependence on $\gamma_{0,T}$]
  Since $B_T=B_0e^{\beta\gamma_1T}$ grows exponentially in $T$, the constant $\gamma_{0,T}$ can be very small for motilities that decay to zero, making the prefactor in \eqref{eq:gronwall-estimate} very large. The constant enters the proof because the mixed term in \eqref{eq:mu-e-splitting} and the residual $R_c$ are paired with $e$ in $L^2(\Omega)$, and the resulting factor $\|e\|_{L^2(\Omega)}$ is absorbed by the coercive term~\eqref{eq:coercive}. Estimating these two terms instead by $\|e\|_{H^{-1}(\Omega)}$ times the $L^2(\Omega)$-norm of the gradient of the remaining factor would avoid this absorption, but it requires $(\gamma(v)-\gamma(\bar v))\bar\rho$ and $R_c$ to belong to $H^1(\Omega)$. The density reconstruction that we consider in Section~\ref{sec:reconstruction} and that allows for mesh adaptivity is constant on each cell and hence does not have this regularity; see Remark~\ref{rem:Rc-regularity} for $R_c$.

  We also note that simulations of interest including our experiments in Section~\ref{sec:experiments} require large $\beta$, cf.~\cite{2025HerdaTrescasesZurek}, which can lead to infeasibly large $B_T$ in practice. In those cases we work with upper bounds $\tilde B_T  \geq v$ that we estimate from numerical experiments and take $\gamma_{0,T} = \min_{0 \leq s \leq \tilde B_T} \gamma(s)$. In this case the stability estimate of Theorem~\ref{thm:hminus1-stability} remains valid as long as the estimated upper bound is not violated.
\end{remark}

\section{Numerical scheme}\label{sec:scheme}
This section discusses the finite volume scheme that we employ in this work.
It is a quasi-stationary adaptation of the fully discrete scheme proposed in
\cite{2025HerdaTrescasesZurek}, which uses a two-point flux approximation
in space and a semi-implicit Euler method in time.

Let $\mathcal T_h$ be an admissible finite volume mesh over
$\Omega$, consisting of cells $K$, faces $\sigma$, and associated cell
midpoints $x_K$, cf.~\cite{eymard2000finitvolummethod}. Let $\mathcal E_h$ denote the set of mesh faces forming the union
\[
\mathcal E_h
=
\mathcal E_h^{\mathrm{int}}
\mathbin{\dot\cup}
\mathcal E_h^{\mathrm{ext}},
\]
where $\mathcal E_h^{\mathrm{int}}$ consists of the interior faces
$\sigma=K|L$ separating two cells $K, L \in \mathcal T_h$, while $\mathcal E_h^{\mathrm{ext}}$ consists of the faces
contained in $\partial\Omega$. For each $K\in\mathcal T_h$, let
$\mathcal E_K$ denote the set of its faces and further define
\[
\mathcal E_K^{\mathrm{int}}
\defeq
\mathcal E_K\cap\mathcal E_h^{\mathrm{int}},
\qquad
\mathcal E_K^{\mathrm{ext}}
\defeq
\mathcal E_K\cap\mathcal E_h^{\mathrm{ext}}.
\]
For every interior face $\sigma=K|L$, the segment $\overline{x_Kx_L}$ is orthogonal
to $\sigma$. We define the \emph{transmissibility} $\tau_\sigma$ of the face $\sigma$ as 
\[
\tau_\sigma\defeq\frac{|\sigma|}{d_\sigma},  \qquad
d_\sigma
\defeq
\begin{cases}
|x_L-x_K|,
& \sigma=K|L\in\mathcal E_h^{\mathrm{int}},\\[2mm]
\operatorname{dist}(x_K,\sigma),
& \sigma\in\mathcal E_K^{\mathrm{ext}}
\end{cases}.
\]
In dimension $d=1$, we use the convention $|\sigma|=1$. We assume the mesh regularity condition of \cite{2025HerdaTrescasesZurek}: there exists $\zeta>0$, independent of $h$, such that
\[
\operatorname{dist}(x_K,\sigma) \ge \zeta d_\sigma
\qquad \text{for all } K \in \mathcal T_h,\quad \sigma \in \mathcal E_K.
\]
For a cellwise
constant quantity $w_h=(w_K)_{K\in\mathcal T_h}$, we define the difference
\[
D_{K,\sigma}w_h
\defeq
\begin{cases}
w_L-w_K,
& \sigma=K|L\in\mathcal E_K^{\mathrm{int}},\\[2mm]
0,
& \sigma\in\mathcal E_K^{\mathrm{ext}}.
\end{cases}
\]
With this convention, every two-point numerical flux through a boundary
face is zero, which imposes the homogeneous Neumann boundary condition.

Regarding the time discretization, let
\[
0=t^0<t^1<\cdots<t^N=T
\]
be a partition of the interval $[0,T]$. We introduce the notation
\[
\Delta t^n \defeq t^{n+1}-t^n,
\qquad
I_n \defeq (t^n,t^{n+1}]
\qquad
\text{for } n=0,\dots,N-1.
\]
For $n=0,\dots,N$, we
denote the approximate volume average of the cell density over $K$ at time $t^n$ by $\rho_K^n$, and similarly
we denote the discrete chemoattractant value in $K$ by $v_K^n$.
The discrete initial cell density is given by the cell averages $\rho_K^0\defeq\frac{1}{|K|}\int_K\rho^0\,dx$. We write
$\rho_h^n$ and $v_h^n$ for the associated cellwise constant functions
on $\Omega$.

\paragraph{Solution update.} At each time instant $t^n$ for $n=0,\dots,N$, we first define $v_h^n$ as the solution of the linear system
\begin{equation}\label{eq:scheme-elliptic}
-\sum_{\sigma\in\mathcal E_K^{\mathrm{int}}}
\tau_\sigma D_{K,\sigma}v_h^n
+ \beta |K|v_K^n
= \alpha |K|\rho_K^n
\qquad \text{for all } K \in \mathcal T_h.
\end{equation}
The cellwise flux potential
\[
m_K^{n+1/2} \defeq \rho_K^{n+1}\gamma(v_K^n)
\]
is implicit in the density and explicit in the chemoattractant.
For each $\sigma\in\mathcal E_K$, the numerical flux is given by
\begin{equation}\label{eq:face-flux}
F_{K,\sigma}^{n+1/2}
\defeq
\begin{cases}
-\tau_\sigma D_{K,\sigma}m^{n+1/2}
=-\tau_\sigma(m_L^{n+1/2}-m_K^{n+1/2}),
& \sigma=K|L\in\mathcal E_K^{\mathrm{int}},\\[2mm]
0,
& \sigma\in\mathcal E_K^{\mathrm{ext}}.
\end{cases}
\end{equation}
The cell density update is finally determined by the discrete local conservation law
\begin{equation}\label{eq:scheme-density}
|K|\frac{\rho_K^{n+1}-\rho_K^n}{\Delta t^n}
+\sum_{\sigma\in\mathcal E_K}F_{K,\sigma}^{n+1/2}
=0
\qquad \text{for all } K\in\mathcal T_h,
\end{equation}
which gives rise to a linear system in $\rho_h^{n+1}$.

In~\cite{2025HerdaTrescasesZurek}, the \emph{positivity-preserving property} of this scheme was shown: if the initial cell averages are nonnegative, then
\[
\rho_K^n \ge 0
\qquad \text{for all } K \in \mathcal T_h,\quad n=0,\dots,N.
\]
The discrete chemoattractant also remains nonnegative, as we can easily show. Since boundary faces do not contribute to \eqref{eq:scheme-elliptic}, the discrete elliptic equation can be written as
\[
\left(
\beta |K|
+\sum_{\sigma\in\mathcal E_K^{\mathrm{int}}}\tau_\sigma
\right)v_K^n
-
\sum_{\sigma=K|L\in\mathcal E_K^{\mathrm{int}}}
\tau_\sigma v_L^n
= \alpha |K|\rho_K^n.
\]
The coefficient matrix $A=(A_{KL})$ of this system has positive diagonal entries and nonpositive off-diagonal entries. It is strictly diagonally dominant because
\[
A_{KK}-\sum_{L\ne K}|A_{KL}|
=
\beta |K|>0.
\]
Therefore, $A$ is a nonsingular $M$-matrix. Since $\alpha|K|\rho_K^n\ge0$,
the inverse positivity of $A$ implies
\[
v_K^n\ge0
\qquad
\text{for all }K\in\mathcal T_h,\quad n=0,\dots,N.
\]
Moreover, in a cell $K$ in which $v_h^n$ attains its maximum, the discrete elliptic equation gives $\beta|K|v_K^n\le\alpha|K|\rho_K^n$. Hence
\begin{equation}\label{eq:discrete-signal-bound}
\max_{K\in\mathcal T_h}v_K^n
\le
\frac{\alpha}{\beta}\max_{K\in\mathcal T_h}\rho_K^n
\qquad
\text{for }n=0,\dots,N.
\end{equation}

\section{Reconstruction}\label{sec:reconstruction}
In this section, we introduce reconstructions of the finite volume solution
computed with the scheme in Section~\ref{sec:scheme}.
We apply the stability estimate in Theorem~\ref{thm:hminus1-stability} to the
resulting approximation to derive a posteriori error estimates for the scheme. The stability framework requires a density approximation $\bar\rho$ and a flux approximation $\bar J$ that satisfy the continuity equation \eqref{eq:recon-mass}, a signal approximation $\bar v$ that satisfies the elliptic equation \eqref{eq:recon-elliptic}, and a flux potential approximation $\bar m$ with values in $H^1(\Omega)$. We obtain $\bar\rho$ by affine interpolation of the cell values in time, $\bar J$ from the numerical fluxes through local mixed problems, $\bar m$ as a conforming interpolation of the cellwise flux potential, and $\bar v$ as the elliptic reconstruction associated with $\bar\rho$. The potential and flux reconstructions are  motivated by a posteriori estimates for the heat equation in \cite[Section~3.1.1]{2010ErnVohralik}.

\subsection{\texorpdfstring{Density reconstruction and discrete chemoattractant interpolation}{Density reconstruction and discrete chemoattractant interpolation}}\label{subsec:density-reconstruction}

We start by discussing the reconstruction of the cell densities. For each interval $I_n$ and each cell $K\in\mathcal T_h$, we define the affine density reconstruction
\[
\theta_n(t) \defeq \frac{t-t^n}{\Delta t^n},
\qquad
\bar\rho(t)\big|_K
\defeq
(1-\theta_n(t))\rho_K^n+\theta_n(t)\rho_K^{n+1}
\qquad \text{for } t\in I_n,
\]
and at the initial time we set
\[
\bar\rho(0)\big|_K \defeq \rho_K^0.
\]
Then $\bar\rho$ is piecewise affine in time and cellwise constant in space, and its time derivative is given by
\[
\partial_t\bar\rho(t)\big|_K = \frac{\rho_K^{n+1}-\rho_K^n}{\Delta t^n}
\qquad \text{for } t \in I_n.
\]
Using this reconstruction, the finite volume scheme \eqref{eq:scheme-density} can be written as
\[
|K|\,\partial_t\bar\rho(t)\big|_K + \sum_{\sigma \in \mathcal E_K}F_{K,\sigma}^{n+1/2} = 0
\qquad \text{for } t \in I_n,
\]
which is the local conservation law associated with the numerical fluxes \eqref{eq:face-flux}.

We next define the affine-in-time interpolation of the discrete chemoattractant
\begin{equation}\label{eq:vtilde-def}
\widetilde v_h(t)\big|_K \defeq (1-\theta_n(t))v_K^n + \theta_n(t)v_K^{n+1}
\qquad \text{for } t \in I_n,
\end{equation}
and at the initial time we set
\[
\widetilde v_h(0)\big|_K\defeq v_K^0.
\]

\subsection{\texorpdfstring{Flux potential reconstruction}{Flux potential reconstruction}}\label{subsec:flux-potential-reconstruction}

\begin{figure}
  \centering
  \iftoggle{usetikz}{\input{figs/conforming-flux-potential-reconstruction.tikz}}{\includegraphics{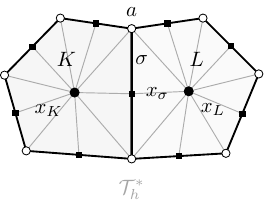}}
    \hspace{8em}
    \iftoggle{usetikz}{\input{figs/tensor-refinement.tikz}}{\includegraphics{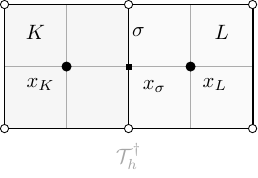}}
  \caption{Left: simplicial refinement $\mathcal T_h^\ast$ for the conforming $\mathbb P_1$ flux potential reconstruction $\bar m^{n+1/2}$. The nodal values are prescribed at the cell centers, face barycenters, and mesh vertices, yielding a globally continuous piecewise affine reconstruction on $\mathcal T_h^\ast$. Right: tensor-product refinement on a Cartesian mesh allowing for a conforming $\mathbb Q_1$ flux potential reconstruction $\bar m^{n+1/2}$ in Remark~\ref{rem:tensor-product-auxiliary-meshes}.}\label{fig:conforming-flux-potential-reconstruction}
\end{figure}

To control the flux residual $R_f=\bar J+\nabla\bar m$, we first reconstruct the cellwise flux potential $m^{n+1/2}$ as a globally continuous function $\bar m^{n+1/2}$. For this purpose, we choose a conforming simplicial refinement $\mathcal T_h^\ast$ subordinate to the finite volume mesh $\mathcal T_h$. This auxiliary mesh is obtained by subdividing each cell into simplices whose nodes include the cell centers $x_K$, the face barycenters $x_\sigma$, and the mesh vertices $a$. It allows us to prescribe single-valued nodal values and thereby obtain a conforming piecewise affine reconstruction. Fig.~\ref{fig:conforming-flux-potential-reconstruction} (left) illustrates this local construction.

We define the flux potential reconstruction $\bar m^{n+1/2}$ by
assigning its values at the cell centers and face barycenters as follows:
\[
\bar m^{n+1/2}(x_K)=m_K^{n+1/2},
\]
\[
\bar m^{n+1/2}(x_\sigma)
=
\begin{cases}
\dfrac{m_K^{n+1/2}+m_L^{n+1/2}}{2},
& \sigma=K|L\in\mathcal E_h^{\mathrm{int}},\\[3mm]
m_K^{n+1/2},
& \sigma\in\mathcal E_K^{\mathrm{ext}}.
\end{cases}
\]
The boundary-face value is consistent with the convention
\[
D_{K,\sigma}m^{n+1/2}=0
\qquad
\text{for }\sigma\in\mathcal E_K^{\mathrm{ext}}.
\]
For each mesh vertex $a$, let
\[
\omega_a
\defeq
\left\{
K\in\mathcal T_h
\mid
a\in\overline K
\right\},
\]
and set
\[
\bar m^{n+1/2}(a)
=
\frac{1}{\#\omega_a}
\sum_{K\in\omega_a}m_K^{n+1/2}.
\]
These single-valued nodal data determine a unique globally continuous
piecewise affine function $\bar m^{n+1/2}$ on $\mathcal T_h^\ast$. Hence, it holds
\[
\bar m^{n+1/2} \in H^1(\Omega),
\qquad
\bar m^{n+1/2}|_S \in \mathbb P_1(S)
\quad \text{for every } S \in \mathcal T_h^\ast.
\]

\begin{remark}[Tensor-product auxiliary meshes]\label{rem:tensor-product-auxiliary-meshes}
  The numerical evaluation of the final error estimates is simpler and less computationally expensive on Cartesian meshes. On such meshes, we define $\bar m^{n+1/2}$ as the conforming $\mathbb Q_1$ function on the tensor-product auxiliary refinement that we denote $\mathcal{T}_h^\dag$. This reduces the number of subcells compared with the simplicial refinement. In this construction, illustrated in Fig.~\ref{fig:conforming-flux-potential-reconstruction} (right), the value of the reconstruction at each node is the average of $m_K^{n+1/2}$ over the cells $K$ whose closure contains the node. This flux potential reconstruction also belongs to $H^1(\Omega)$, which is the only property of $\bar m^{n+1/2}$ used in Theorem~\ref{thm:semi-computable-estimate}. Because the scheme requires the stated orthogonality assumption, it cannot be used with quadtree-/octree-based refinement on such meshes in dimensions $d\in\{2,3\}$; developing an adaptive scheme for this purpose is beyond the scope of this work.
\end{remark}

\subsection{\texorpdfstring{Flux reconstruction}{Flux reconstruction}}\label{subsec:flux-reconstruction}
We construct the flux reconstruction $\bar J^{n+1/2}$ so that the discrete
local conservation law yields
\[
\partial_t\bar\rho+\nabla\cdot\bar J=0.
\]
The locally conservative $H(\mathrm{div})$ flux reconstruction based on local mixed
problems that we introduce in the following is motivated by similar reconstructions used in a posteriori estimates for the heat equation in 
\cite[Section~4.4]{2010ErnVohralik}.
For each $K\in\mathcal T_h$, let
\begin{equation}\label{eq:local-submesh}
\mathcal T_h^\ast(K)
\defeq
\left\{
S\in\mathcal T_h^\ast
\mid
S\subset K
\right\}
\end{equation}
be the induced local simplicial submesh. We denote the associated lowest-order
Raviart--Thomas and piecewise constant spaces by
$RT_0(\mathcal T_h^\ast(K))$ and $P_0(\mathcal T_h^\ast(K))$, respectively,
and set
\[
RT_{0,0}(\mathcal T_h^\ast(K))
\defeq
\left\{
w_h\in RT_0(\mathcal T_h^\ast(K))
\,\middle|\,
w_h\cdot \nu_K=0
\text{ on }\partial K
\right\}.
\]
In the next step, we define
\[
\left.\bar J^{n+1/2}\right|_K
\defeq
q_{K,h}^{n+1/2},
\]
where the pair
\[
(q_{K,h}^{n+1/2},p_{K,h}^{n+1/2})
\in
RT_0(\mathcal T_h^\ast(K))
\times
P_0(\mathcal T_h^\ast(K))
\]
is determined by the following local mixed problem. The normal component of $q_{K,h}^{n+1/2}$ on the boundary faces of $\mathcal T_h^\ast(K)$ is prescribed by the numerical fluxes,
\[
\left.q_{K,h}^{n+1/2}\cdot \nu_K\right|_\sigma
=
\begin{cases}
\dfrac{F_{K,\sigma}^{n+1/2}}{|\sigma|},
&
\sigma\in\mathcal E_K^{\mathrm{int}},
\\[2mm]
0,
&
\sigma\in\mathcal E_K^{\mathrm{ext}},
\end{cases}
\]
and the pair satisfies
\[
\int_K q_{K,h}^{n+1/2}\cdot\phi_h\,dx
-
\int_K p_{K,h}^{n+1/2}\nabla\cdot\phi_h\,dx
=0
\qquad
\text{for every }\phi_h\in RT_{0,0}(\mathcal T_h^\ast(K)),
\]
and
\[
\int_K(\nabla\cdot q_{K,h}^{n+1/2})\psi_h\,dx
=
-\int_K
\left.\partial_t\bar\rho\right|_{K\times I_n}
\psi_h\,dx
\qquad
\text{for every }\psi_h\in P_0(\mathcal T_h^\ast(K)).
\]
Since the numerical flux vanishes on $\mathcal E_K^{\mathrm{ext}}$, the
finite volume scheme implies
\[
\sum_{\sigma\in\mathcal E_K}
F_{K,\sigma}^{n+1/2}
=
-|K|
\left.\partial_t\bar\rho\right|_{K\times I_n}.
\]
This is precisely the compatibility condition for the prescribed normal
trace and divergence. Thus the local mixed problem admits a solution, for which $q_{K,h}^{n+1/2}$ is unique, whereas $p_{K,h}^{n+1/2}$ is
determined up to an additive constant. In particular, it holds
\[
\nabla\cdot\bar J^{n+1/2}
=
-\left.\partial_t\bar\rho\right|_{K\times I_n}
\quad\text{in }K
\quad\text{and} \quad
\bar J^{n+1/2}\cdot \nu_K
=
0
\quad\text{on every }
\sigma\in\mathcal E_K^{\mathrm{ext}}.
\]
Next, on each interior face $\sigma=K|L$, we fix an orientation
\[
\nu_\sigma
=
\nu_{K,\sigma}
=
-\nu_{L,\sigma}
\]
Due to the antisymmetry
$F_{K,\sigma}^{n+1/2}
=
-F_{L,\sigma}^{n+1/2}$
and the prescribed normal component, we have
\[
\left.\bar J^{n+1/2}\right|_K\cdot \nu_\sigma
=
\frac{F_{K,\sigma}^{n+1/2}}{|\sigma|}
=
-\frac{F_{L,\sigma}^{n+1/2}}{|\sigma|}
=
\left.\bar J^{n+1/2}\right|_L\cdot \nu_\sigma.
\]
Hence, the normal component of $\bar J^{n+1/2}$ agrees across every interior
face. For every boundary face $\sigma\in\mathcal E_K^{\mathrm{ext}}$, we
have $\nu_K=\nu$, and the boundary flux condition implies
\[
\left.\bar J^{n+1/2}\right|_K\cdot\nu
=
\frac{F_{K,\sigma}^{n+1/2}}{|\sigma|}
=
0.
\]
Because $\left.\bar J^{n+1/2}\right|_K\in H(\mathrm{div};K)$ for every
$K\in\mathcal T_h$, the continuity across interior faces and the zero normal
component on $\partial\Omega$ imply
\[
\bar J^{n+1/2}
\in
H_0(\mathrm{div};\Omega),
\]
where
\[
H_0(\mathrm{div};\Omega)
\defeq
\left\{
q\in H(\mathrm{div};\Omega)
\,\middle|\,
q\cdot\nu=0
\text{ on }\partial\Omega
\right\}.
\]
For $t\in I_n$, we extend the flux potential reconstruction and the flux
reconstruction in time by taking
\[
\bar m(t,\cdot)\defeq\bar m^{n+1/2},
\qquad
\bar J(t,\cdot)\defeq\bar J^{n+1/2}.
\]
\begin{remark}
The reconstructions satisfy
\[
\bar m\in L^\infty(0,T;H^1(\Omega)),
\qquad
\bar J\in L^2(0,T;H_0(\mathrm{div};\Omega))
\]
and the local divergence identities yield the global continuity equation
\[
\partial_t\bar\rho+\nabla\cdot\bar J=0
\qquad
\text{in }L^2(0,T;L^2(\Omega)).
\]
\end{remark}

\begin{remark}[Explicit flux reconstructions]\label{rem:cartesian-flux-reconstruction}
  \begin{enumerate}
    \item Given a mesh cell $K=(a,b)$ in one space dimension the above construction implies the explicit form of the flux reconstruction
\[
  \left.\bar J^{n+1/2}\right|_K = -F_{K, \sigma_a}^{n+1/2} + \frac{F_{K, \sigma_a}^{n+1/2} + F_{K, \sigma_b}^{n+1/2}}{b-a} \, (x-a) \quad \text{for }x \in \bar K,
\]
where $\sigma_a$ and $\sigma_b$ denote the left and the right faces of $K$. 
\item On Cartesian meshes in two or three space dimensions using tensor-product refinement (see Remark~\ref{rem:tensor-product-auxiliary-meshes}) we employ local submeshes $\mathcal T_h^\dag(K)$ in analogy to \eqref{eq:local-submesh}. The solution to the local mixed problem in $RT_{0,0}(\mathcal T_h^\dag(K))$ can then also be explicitly determined. For $K=(a,b)\times(c,d)$ in two space dimensions, we obtain
\[
  \left.\bar J^{n+1/2}\right|_K =
  \begin{pmatrix}
    \dfrac{1}{d-c}
    \left(
    -F_{K, \sigma_W}^{n+1/2} + \dfrac{F_{K, \sigma_W}^{n+1/2} + F_{K, \sigma_E}^{n+1/2}}{b-a} \, (x-a)
    \right) \\[9pt]
    \dfrac{1}{b-a}
    \left(
    -F_{K, \sigma_S}^{n+1/2} + \dfrac{F_{K, \sigma_S}^{n+1/2} + F_{K, \sigma_N}^{n+1/2}}{d-c} \, (y-c)
    \right)
  \end{pmatrix}
  \quad \text{for }(x, y) \in \bar K,
\]
where $\sigma_W$, $\sigma_N$, $\sigma_E$ and $\sigma_S$ denote the west, north, east and south faces of $K$. The properties of $\bar J^{n+1/2}$ established in this subsection remain valid.
\end{enumerate}
\end{remark}

\subsection{\texorpdfstring{Elliptic reconstruction}{Elliptic reconstruction}}\label{subsec:elliptic-reconstruction}

We use the elliptic reconstruction introduced in \cite{2003MakridakisNochetto} to
associate the density reconstruction $\bar\rho(t)$ with a chemoattractant
reconstruction $\bar v(t)$. For every $t\in[0,T]$,
$\bar v(t)\in H^1(\Omega)$ satisfies
\[
\int_\Omega \nabla\bar v(t)\cdot\nabla\eta\,dx
+
\beta\int_\Omega \bar v(t)\eta\,dx
=
\alpha\int_\Omega \bar\rho(t)\eta\,dx
\qquad
\text{for every }\eta\in H^1(\Omega).
\]
This weak formulation incorporates the homogeneous Neumann boundary condition
$\partial_\nu\bar v(t)=0$ on $\partial\Omega$. Since $\beta>0$, there exists
a unique $\bar v\in L^\infty(0,T;H^1(\Omega))$ satisfying the weak
formulation.

The density reconstruction is nonnegative because it interpolates nonnegative
cell values, and it satisfies
\[
\|\bar\rho\|_{L^\infty((0,T)\times\Omega)}
=
\max_{0\le n\le N}\max_{K\in\mathcal T_h}\rho_K^n.
\]
Since the constant $\alpha\|\bar\rho(t)\|_{L^\infty(\Omega)}/\beta$ solves the
same weak problem with $\bar\rho(t)$ replaced by
$\|\bar\rho(t)\|_{L^\infty(\Omega)}$, the weak maximum principle gives
\[
0\le\bar v(t)\le\frac{\alpha}{\beta}\|\bar\rho(t)\|_{L^\infty(\Omega)}
\qquad
\text{almost everywhere in }\Omega.
\]
The interpolation $\widetilde v_h$ defined in \eqref{eq:vtilde-def} is also nonnegative, because all discrete
chemoattractant states are nonnegative, and the discrete bound
\eqref{eq:discrete-signal-bound} gives
\[
\widetilde v_h(t)\le\frac{\alpha}{\beta}\|\bar\rho\|_{L^\infty((0,T)\times\Omega)}
\qquad
\text{for every }t\in[0,T].
\]
Hence $\gamma(\bar v)$ and
$\gamma(\widetilde v_h)$ are well-defined under
\eqref{eq:gamma-assumptions}. The regularity of $\bar m$, $\bar J$, and
$\bar\rho$ yields
\[
\begin{aligned}
R_f
&=
\bar J+\nabla\bar m
\in
L^2(0,T;L^2(\Omega;\R^d)),
\\
R_c
&=
\bar m-\gamma(\bar v)\bar\rho
\in
L^2((0,T)\times\Omega).
\end{aligned}
\]

\begin{remark}[Regularity of $R_c$]\label{rem:Rc-regularity}
For almost every $t\in(0,T)$, the residual $R_c(t)$ belongs to $H^1(K)$ for
every $K\in\mathcal T_h$, because $\bar m(t)$ and $\gamma(\bar v(t))$ belong to
$H^1(\Omega)$ and $\bar\rho(t)$ is constant on $K$. Across an interior face
$\sigma=K|L$, the traces of $R_c(t)$ from $K$ and $L$ differ by
$\gamma(\bar v(t))\bigl(\bar\rho(t)|_L-\bar\rho(t)|_K\bigr)$. Since $\gamma>0$,
$R_c(t)$ therefore belongs neither to $H^1(\Omega)$ nor to $W^{1,1}(\Omega)$
unless $\bar\rho(t)$ is constant on $\Omega$.
\end{remark}

\begin{remark}[Adaptivity]\label{rem:reconstruction-adaptivity}
  Our choice of reconstruction is stable under mesh refinement in the sense that a projection of $\rho_h$ to a locally refined mesh does not change the reconstruction~$\bar \rho$. This property has enabled adaptive mesh refinement that also preserves the total mass; see Section~\ref{sec:refinementexperiment}. In \cite{2025GiesselmannKolbe}, where the cell density was reconstructed using linear interpolation on the dual mesh, an analogous adaptive refinement strategy was not feasible because such a stability property was unavailable.

  While coarsening mesh cells does not preserve the reconstruction in our approach, this has not caused issues with the adaptive strategy in practice. The estimator-based mesh update we employ has ensured that only regions with very small observed gradients of $\rho$ were coarsened, leading to only negligible changes in $\bar \rho$ after projection. 
\end{remark}

\section{A posteriori error estimate}
In this section, we conclude the analytical part of this work and establish an a posteriori error estimate for the scheme from Section~\ref{sec:scheme}.

\subsection{\texorpdfstring{Estimator for the elliptic reconstruction error}{Estimator for the elliptic reconstruction error}}\label{subsec:auxiliary-elliptic-control}

The elliptic reconstruction $\bar v$ is not available explicitly, and the
discrete chemoattractant $\widetilde v_h$ is cellwise constant, so residual
estimates for conforming finite elements do not apply to
$\bar v-\widetilde v_h$ directly. To estimate
$\|\bar v-\widetilde v_h\|_{L^2(\Omega)}$, we therefore introduce a conforming
Galerkin approximation $\widetilde v_{h,a}$ of $\bar v$ on the simplicial
refinement $\mathcal T_h^\ast$. We bound
$\|\bar v-\widetilde v_{h,a}\|_{L^2(\Omega)}$ by a residual estimator and
compute $\|\widetilde v_{h,a}-\widetilde v_h\|_{L^2(\Omega)}$ directly. Let
\[
V_h^a
\defeq
\left\{
\phi_h\in H^1(\Omega)
\,:\,
\left.\phi_h\right|_S\in\mathbb P_1(S)
\text{ for every }S\in\mathcal T_h^\ast
\right\}.
\]
The homogeneous Neumann condition is incorporated in the weak problem, so no
boundary constraint is imposed on $V_h^a$. For each $n=0,\dots,N$, let
$v_{h,a}^n\in V_h^a$ be the unique function satisfying
\begin{equation}\label{eq:aux-conforming-elliptic}
\int_\Omega
\nabla v_{h,a}^n\cdot\nabla\phi_h\,dx
+
\beta\int_\Omega
v_{h,a}^n\phi_h\,dx
=
\alpha\int_\Omega
\rho_h^n\phi_h\,dx
\qquad
\text{for every }\phi_h\in V_h^a.
\end{equation}
For $t\in I_n$, define
\[
\widetilde v_{h,a}(t)
\defeq
(1-\theta_n(t))v_{h,a}^n
+
\theta_n(t)v_{h,a}^{n+1},
\]
and set $\widetilde v_{h,a}(0)\defeq v_{h,a}^0$. The definition of
$\bar\rho$ and the linearity of the discrete weak problem give
\begin{equation}\label{eq:aux-conforming-elliptic-time}
\int_\Omega
\nabla\widetilde v_{h,a}(t)\cdot\nabla\phi_h\,dx
+
\beta\int_\Omega
\widetilde v_{h,a}(t)\phi_h\,dx
=
\alpha\int_\Omega
\bar\rho(t)\phi_h\,dx
\qquad
\text{for every }\phi_h\in V_h^a.
\end{equation}
Thus $\widetilde v_{h,a}(t)$ is the conforming Galerkin approximation of the
elliptic reconstruction $\bar v(t)$.

Let $\mathcal F_h^{\mathrm{int}}$ and $\mathcal F_h^\partial$ denote the sets of interior and boundary faces of the simplicial refinement $\mathcal T_h^\ast$, respectively. For $S\in\mathcal T_h^\ast$, set
\[
h_S \defeq \operatorname{diam}(S),
\]
and define
\[
h_F
\defeq
\begin{cases}
\min\{h_{S^+},h_{S^-}\},
&
F=\partial S^+\cap\partial S^-
\in\mathcal F_h^{\mathrm{int}},
\\[2mm]
h_S,
&
F\subset\partial S\cap\partial\Omega
\text{ and }F\in\mathcal F_h^\partial.
\end{cases}
\]
For $d=1$, the $L^2(F)$-norm of a quantity defined on a face is interpreted as its absolute value. Choose a unit normal $\nu_F$ on every $F\in\mathcal F_h^{\mathrm{int}}$ and label the adjacent elements $S^+$ and $S^-$ so that $\nu_F$ points from $S^+$ to $S^-$. For a function $w$ that is affine on every $S\in\mathcal T_h^\ast$, we define the jump of its normal derivative across $F$ by
\[
\left.[\nabla w\cdot \nu_F]\right|_F
\defeq
\left.\nabla w\right|_{S^+}\cdot \nu_F
-
\left.\nabla w\right|_{S^-}\cdot \nu_F,
\qquad
F\in\mathcal F_h^{\mathrm{int}}.
\]
For $F\in\mathcal F_h^\partial$, the normal derivative
\[
\left.\nabla w\cdot\nu\right|_F,
\]
where $\nu$ is the outward unit normal to $\partial\Omega$, measures the violation of the homogeneous Neumann condition. The orientation selected on an interior face does not affect the estimator in Proposition~\ref{prop:auxiliary-elliptic-control} because the jump appears only through its $L^2(F)$-norm.

\begin{proposition}[Elliptic reconstruction error estimator]\label{prop:auxiliary-elliptic-control}
For almost every $t\in(0,T)$, define
\begin{equation}\label{eq:explicit-elliptic-indicator}
\begin{aligned}
\eta_0(t)^2
\defeq{}&
\sum_{S\in\mathcal T_h^\ast}
h_S^4
\|\alpha\bar\rho(t)-\beta\widetilde v_{h,a}(t)\|_{L^2(S)}^2
\\
&+
\sum_{F\in\mathcal F_h^{\mathrm{int}}}
h_F^3
\|[\nabla\widetilde v_{h,a}(t)\cdot \nu_F]\|_{L^2(F)}^2
\\
&+
\sum_{F\in\mathcal F_h^\partial}
h_F^3
\|\nabla\widetilde v_{h,a}(t)\cdot\nu\|_{L^2(F)}^2.
\end{aligned}
\end{equation}
Then there exists a constant $C_0>0$, depending only on $\beta$, $\Omega$
and the shape regularity of $\mathcal T_h^\ast$, such that
\begin{equation}\label{eq:explicit-elliptic-estimator}
\|\bar v(t)-\widetilde v_{h,a}(t)\|_{L^2(\Omega)}
\le
C_0\,\eta_0(t)
\qquad
\text{for almost every }t\in(0,T).
\end{equation}
\end{proposition}
\begin{proof}
Fix $t\in I_n$ and set $e_v\defeq\bar v(t)-\widetilde v_{h,a}(t)$. Subtracting
\eqref{eq:aux-conforming-elliptic-time} from the weak problem for $\bar v(t)$
gives the Galerkin orthogonality
\[
\int_\Omega\nabla e_v\cdot\nabla\phi_h\,dx
+\beta\int_\Omega e_v\phi_h\,dx
=0
\qquad\text{for every }\phi_h\in V_h^a.
\]
We use the duality argument for residual estimates in the $L^2$-norm; see
\cite[Section~4]{2003MakridakisNochetto} and, for the Poisson problem,
\cite[Section~1.11]{2013Verfurth}. Let $z\in H^1(\Omega)$ be the unique
solution of
\[
\int_\Omega\nabla z\cdot\nabla\psi\,dx
+\beta\int_\Omega z\psi\,dx
=\int_\Omega e_v\psi\,dx
\qquad\text{for every }\psi\in H^1(\Omega).
\]
Since $-\Delta z=e_v-\beta z$ has zero mean and $\Omega$ is convex, homogeneous
Neumann regularity \cite[Theorem~2.4.3 and Remark~2.6.9]{1992Grisvard} and the
energy estimate give a constant $C_{\mathrm{ell}}>0$, depending only on $\beta$
and $\Omega$, such that
\[
\|z\|_{H^2(\Omega)}\le C_{\mathrm{ell}}\|e_v\|_{L^2(\Omega)}.
\]
Let $I_h$ denote the Scott--Zhang interpolation operator onto $V_h^a$
\cite{1990ScottZhang} and set $w\defeq z-I_hz$. Choosing $\psi=e_v$ in the dual
problem, using the Galerkin orthogonality, and integrating by parts on each
$S\in\mathcal T_h^\ast$, where $\Delta\widetilde v_{h,a}(t)=0$, we obtain
\[
\begin{aligned}
\|e_v\|_{L^2(\Omega)}^2
={}&
\sum_{S\in\mathcal T_h^\ast}\int_S
\bigl(\alpha\bar\rho(t)-\beta\widetilde v_{h,a}(t)\bigr)w\,dx
-\sum_{F\in\mathcal F_h^{\mathrm{int}}}\int_F
[\nabla\widetilde v_{h,a}(t)\cdot\nu_F]\,w\,ds
\\
&-\sum_{F\in\mathcal F_h^\partial}\int_F
\bigl(\nabla\widetilde v_{h,a}(t)\cdot\nu\bigr)w\,ds.
\end{aligned}
\]
Let $\omega_S$ denote the union of the elements sharing a vertex with $S$, and
let $\omega_F$ denote the union of $\omega_S$ over the elements $S$ with
$F\subset\partial S$. The interpolation estimates for $I_h$ and the scaled
trace inequality give a constant $C_I>0$, depending only on $d$ and the shape
regularity of $\mathcal T_h^\ast$, such that
\[
\|w\|_{L^2(S)}\le C_I h_S^2|z|_{H^2(\omega_S)},
\qquad
\|w\|_{L^2(F)}\le C_I h_F^{3/2}|z|_{H^2(\omega_F)}
\]
for every element $S$ and every face $F$. Since the patches overlap a bounded
number of times, the sum of $|z|_{H^2(\omega_S)}^2$ and
$|z|_{H^2(\omega_F)}^2$ over all elements and faces is bounded by
$C_{\mathrm{ov}}\|z\|_{H^2(\Omega)}^2$, where $C_{\mathrm{ov}}>0$ depends only
on $d$ and the shape regularity of $\mathcal T_h^\ast$. The Cauchy--Schwarz
inequality and the definition of $\eta_0(t)$ therefore give
\[
\|e_v\|_{L^2(\Omega)}^2
\le C_I C_{\mathrm{ov}}^{1/2}\,\eta_0(t)\,\|z\|_{H^2(\Omega)}
\le C_I C_{\mathrm{ov}}^{1/2}C_{\mathrm{ell}}\,\eta_0(t)\,\|e_v\|_{L^2(\Omega)},
\]
which is \eqref{eq:explicit-elliptic-estimator} with
$C_0\defeq C_I C_{\mathrm{ov}}^{1/2}C_{\mathrm{ell}}$.
\end{proof}

\begin{remark}[Control of $\bar v-\widetilde v_h$]
Proposition~\ref{prop:auxiliary-elliptic-control} and the triangle inequality imply
\[
\begin{aligned}
\|\bar v(t)-\widetilde v_h(t)\|_{L^2(\Omega)}
\le{}&
C_0\eta_0(t)+
\|\widetilde v_{h,a}(t)-\widetilde v_h(t)\|_{L^2(\Omega)},
\end{aligned}
\]
where the difference
$\widetilde v_{h,a}(t)-\widetilde v_h(t)$ is computable.
\end{remark}

\begin{remark}[Computable constants]
To estimate $C_0=C_I C_{\mathrm{ov}}^{1/2}C_{\mathrm{ell}}$
in Proposition~\ref{prop:auxiliary-elliptic-control}, we need explicit upper bounds
for the interpolation and trace constant $C_I$ and the Neumann regularity
constant $C_{\mathrm{ell}}$ for $-\Delta+\beta$.
The overlap constant $C_{\mathrm{ov}}$ can be bounded by the maximum number
of element and face patches containing any one auxiliary element.
Related techniques for evaluating estimator constants are discussed in
\cite[Remarks~6.5 and~7.3]{2024GiesselmannKwon}
and the references therein.
\end{remark}

\begin{remark}[Cartesian auxiliary meshes]\label{rem:cartesian-auxiliary-elliptic}
On a Cartesian mesh in dimensions two and three that employs tensor-product refinement (see Remark~\ref{rem:tensor-product-auxiliary-meshes}) we use the
standard conforming $\mathbb Q_1$ auxiliary Galerkin approximation on $\mathcal T_h^\dag$. We retain the notation
$\widetilde v_{h,a}$, $\eta_0$ and $C_0$, using the same residual formula
and diameter weights on this mesh. The estimate in
Proposition~\ref{prop:auxiliary-elliptic-control} then holds with $C_0$
chosen for this auxiliary space. The proof of
Proposition~\ref{prop:auxiliary-elliptic-control} applies using standard
$\mathbb Q_1$ interpolation and the identity $\Delta q=0$ for
$q\in\mathbb Q_1(S)$ on each Cartesian element.
\end{remark}

\subsection{\texorpdfstring{A posteriori error estimate}{A posteriori error estimate}}\label{subsec:aposteriori-error-estimate}
We now derive an a posteriori error estimate for the finite volume approximation. The residual $R_f$ is computable from the reconstructions. The residual $R_c$ contains the elliptic reconstruction $\bar v$, which is not available explicitly; we therefore replace $\gamma(\bar v)$ by $\gamma(\widetilde v_h)$ and estimate the difference.
\begin{theorem}[A posteriori error estimate]\label{thm:semi-computable-estimate}
Let $(\rho,v)$ be a weak solution to \eqref{eq:model}, and let
$(\bar\rho,\bar v,\bar m,\bar J)$ be the reconstructions defined in
Section \ref{sec:reconstruction}. Set
\[
e \defeq \rho-\bar\rho,
\qquad
R_f \defeq \bar J+\nabla \bar m,
\qquad
R_c \defeq \bar m-\gamma(\bar v)\bar\rho.
\]
Let $B_T$ and $\gamma_{0,T}$ be the constants in Proposition \ref{prop:signal-bound}, set
\[
M_T \defeq \|\bar\rho\|_{L^\infty((0,T)\times\Omega)},
\qquad
\tilde \alpha \defeq \frac{\alpha}{\min\{\beta, 1\}},
\]
and let $L_{\gamma,T}$ denote the Lipschitz constant of $\gamma$ on the interval $[0,\max\{B_T,\alpha M_T/\beta\}]$. Define
\[
\Lambda_T \defeq 1 + \frac{2L_{\gamma,T}^2 M_T^2 \tilde \alpha^2}{\gamma_{0,T}}.
\]
Let $\widetilde v_h$ be the affine interpolation in time
of the discrete chemoattractant, as defined in
Subsection \ref{subsec:density-reconstruction}, and set
\[
M_n \defeq \|\bar\rho\|_{L^\infty(I_n\times \Omega)}
\qquad
\text{for }n=0,\ldots,N-1.
\]
Then
\begin{multline}\label{eq:semi-computable-estimate}
\sup_{0 \le t \le T}\|e(t)\|_{H^{-1}(\Omega)}^2
+ \int_0^T \|e(t)\|_{L^2(\Omega)}^2\,dt \\
\le \left(1+\frac{1}{\gamma_{0,T}}\right)e^{\Lambda_T T}\Biggl(
\|e(0)\|_{H^{-1}(\Omega)}^2
+ \int_0^T \|R_f(t)\|_{L^2(\Omega)}^2\,dt \\
+ \frac{4}{\gamma_{0,T}}\int_0^T \|\bar m(t)-\gamma(\widetilde v_h(t))\bar\rho(t)\|_{L^2(\Omega)}^2\,dt \\
+ \frac{4L_{\gamma,T}^2}{\gamma_{0,T}}\sum_{n=0}^{N-1} M_n^2 \int_{I_n}\|\bar v(t)-\widetilde v_h(t)\|_{L^2(\Omega)}^2\,dt
\Biggr).
\end{multline}
\end{theorem}

\begin{proof}
The regularity and initial condition in Definition \ref{def:weak-solution}
imply
\[
\rho^0-\bar\rho(0)\in(H^1(\Omega))'.
\]
We initialize the cell densities with the cell averages of $\rho^0$, so
\[
\int_\Omega\bar\rho(0,x)\,dx
=
\sum_{K\in\mathcal T_h}|K|\rho_K^0
=
\int_\Omega\rho^0(x)\,dx.
\]
Hence,
\[
\rho^0-\bar\rho(0)\in H^{-1}_\diamond(\Omega).
\]
We use the reconstruction properties established in Section \ref{sec:reconstruction}
to verify the remaining hypotheses of Theorem \ref{thm:hminus1-stability}.
Estimate \eqref{eq:gronwall-estimate} therefore gives
\[
\begin{aligned}
\sup_{0 \le t \le T}\|e(t)\|_{H^{-1}(\Omega)}^2
+ \int_0^T \|e(t)\|_{L^2(\Omega)}^2\,dt
&\le \left(1+\frac{1}{\gamma_{0,T}}\right)e^{\Lambda_T T}\Biggl(
\|e(0)\|_{H^{-1}(\Omega)}^2 \\
&\qquad
+ \int_0^T\|R_f(t)\|_{L^2(\Omega)}^2\,dt \\
&\qquad
+ \frac{2}{\gamma_{0,T}}\int_0^T\|R_c(t)\|_{L^2(\Omega)}^2\,dt
\Biggr).
\end{aligned}
\]
For $n=0,\ldots,N-1$ and almost every $t\in I_n$, we add and subtract
$\gamma(\widetilde v_h(t))\bar\rho(t)$ in $R_c(t)$. By the bounds in
Subsection \ref{subsec:elliptic-reconstruction}, $\bar v(t)$ and
$\widetilde v_h(t)$ take values in $[0,\alpha M_T/\beta]$. The Lipschitz
continuity of $\gamma$ on
this interval, the definition of $M_n$, and
$(a+b)^2 \le 2a^2+2b^2$ imply
\[
\begin{aligned}
\|R_c(t)\|_{L^2(\Omega)}^2
\le{}&
2\|\bar m(t)-\gamma(\widetilde v_h(t))\bar\rho(t)\|_{L^2(\Omega)}^2
\\
&+
2L_{\gamma,T}^2 M_n^2
\|\bar v(t)-\widetilde v_h(t)\|_{L^2(\Omega)}^2.
\end{aligned}
\]
We integrate over $I_n$, sum over $n=0,\ldots,N-1$, and multiply by
$2/\gamma_{0,T}$ to obtain
\[
\begin{aligned}
\frac{2}{\gamma_{0,T}}\int_0^T\|R_c(t)\|_{L^2(\Omega)}^2\,dt
\le\;&
\frac{4}{\gamma_{0,T}}\int_0^T
\|\bar m(t)-\gamma(\widetilde v_h(t))\bar\rho(t)\|_{L^2(\Omega)}^2\,dt
\\
&+
\frac{4L_{\gamma,T}^2}{\gamma_{0,T}}\sum_{n=0}^{N-1} M_n^2
\int_{I_n}\|\bar v(t)-\widetilde v_h(t)\|_{L^2(\Omega)}^2\,dt.
\end{aligned}
\]
Substituting this bound into \eqref{eq:gronwall-estimate}, we obtain
\eqref{eq:semi-computable-estimate}.
\end{proof}

\begin{corollary}[Computable a posteriori error estimate]\label{cor:computable-estimate}
Under the assumptions of Theorem \ref{thm:semi-computable-estimate}, let
$\eta_0$ be defined by \eqref{eq:explicit-elliptic-indicator}, and let $C_0$
be the constant in Proposition \ref{prop:auxiliary-elliptic-control}. Then
\begin{multline}\label{eq:computable-estimate}
\sup_{0 \le t \le T}\|e(t)\|_{H^{-1}(\Omega)}^2
+ \int_0^T \|e(t)\|_{L^2(\Omega)}^2\,dt \\
\le \left(1+\frac{1}{\gamma_{0,T}}\right)e^{\Lambda_T T}\Biggl(
\|e(0)\|_{H^{-1}(\Omega)}^2
+ \int_0^T \|R_f(t)\|_{L^2(\Omega)}^2\,dt \\
+ \frac{4}{\gamma_{0,T}}\int_0^T \|\bar m(t)-\gamma(\widetilde v_h(t))\bar\rho(t)\|_{L^2(\Omega)}^2\,dt \\
+ \frac{8L_{\gamma,T}^2}{\gamma_{0,T}}\sum_{n=0}^{N-1} M_n^2
\int_{I_n}
\Bigl(
C_0^2\eta_0(t)^2
+
\|\widetilde v_{h,a}(t)-\widetilde v_h(t)\|_{L^2(\Omega)}^2
\Bigr)\,dt
\Biggr).
\end{multline}
\end{corollary}

\begin{proof}
For almost every $t\in(0,T)$, Proposition
\ref{prop:auxiliary-elliptic-control} and the triangle inequality give
\[
\|\bar v(t)-\widetilde v_h(t)\|_{L^2(\Omega)}
\le
C_0\eta_0(t)
+
\|\widetilde v_{h,a}(t)-\widetilde v_h(t)\|_{L^2(\Omega)}.
\]
Therefore, $(a+b)^2 \le 2a^2+2b^2$ yields
\[
\|\bar v(t)-\widetilde v_h(t)\|_{L^2(\Omega)}^2
\le
2C_0^2\eta_0(t)^2
+
2\|\widetilde v_{h,a}(t)-\widetilde v_h(t)\|_{L^2(\Omega)}^2.
\]
Substituting this bound into \eqref{eq:semi-computable-estimate}, we obtain
\eqref{eq:computable-estimate}.
\end{proof}
\begin{remark}
  Since our scheme employs a piecewise constant approximation, we generally expect first order accuracy in space on arbitrary grids. If the right-hand side of \eqref{eq:computable-estimate}, which we will call the \emph{estimator}, belongs to the class $\mathcal O (h^2)$ as $h \to 0$, we say that the estimator \emph{scales optimally}, since the same is expected for the squared error of the scheme. This notion disregards the fact that superconvergence was observed in~\cite{2025HerdaTrescasesZurek} in experiments on uniform grids.
\end{remark}

\section{Numerical experiments}\label{sec:experiments}
In this section, we conduct numerical experiments in 1D and 2D using the scheme discussed in Section~\ref{sec:scheme} and analyze the behavior of the derived a posteriori error estimates. Some of our experiments in 1D use adaptive nonuniform grids, and we employ uniform Cartesian grids in the 2D case for simplicity.

\subsection{Preliminaries}
In this section, we discuss some practical aspects of the computations on Cartesian grids and introduce a decomposition of the estimator.

The evaluation of the error estimate in Theorem~\ref{thm:semi-computable-estimate} requires a computation of the reconstruction variables. In all experiments, we use the flux reconstruction of Remark~\ref{rem:cartesian-flux-reconstruction}.

In the following, let $\mathcal{E} \defeq \mathcal{E}[\rho_h; T]$ denote the right-hand side of \eqref{eq:computable-estimate}, i.e., our computable a posteriori bound. We introduce for any $K\in \mathcal{T}_h$ and $n\in \{0,1, \dots, N-1\}$ the following local contributions in space and time:
\[
  \mathcal{E}^n_{1, K} \defeq \int_{I_n} \|\bar J(t) + \nabla \bar m (t) \|_{L^2(K)}^2\,dt, \quad \mathcal{E}_{2, K}^n \defeq \frac{4}{\gamma_{0,T}}\int_{I_n} \|\bar m(t)-\gamma(\widetilde v_h(t))\bar\rho(t)\|_{L^2(K)}^2\,dt
\]
and
\begin{align*}
\mathcal{E}^n_{3,K}
&\defeq \frac{8L_{\gamma,T}^2}{\gamma_{0,T}} M_n^2
\int_{I_n}
\Bigl(
C_0^2\eta_{0,K}(t)^2 + \|\widetilde v_{h,a}(t)-\widetilde v_h(t)\|^2_{L^2(K)}
\Bigr)\,dt
\end{align*}
with the localized elliptic control term
\begin{align*}
\eta_{0,K}(t)^2
&\defeq
\sum_{S \in \mathcal T_h^\ast(K)} h_S^4
\|\alpha\bar\rho(t)-\beta\widetilde v_{h,a}(t)\|_{L^2(S)}^2
+
\sum_{\substack{F \in \mathcal F_h^{\mathrm{int}}\\ F \subset \operatorname{int}(K)}} h_F^3
  \|[\nabla \widetilde v_{h,a}(t)\cdot \nu_F]\|_{L^2(F)}^2 \\
  &\quad + \frac{1}{2} \sum_{\substack{F \in \mathcal F_h^{\mathrm{int}} \\ F \subset \partial K}} h_F^3  \|[\nabla \widetilde v_{h,a}(t)\cdot \nu_F]\|_{L^2(F)}^2 + \sum_{\substack{F \in \mathcal F_h^{\partial} \\ F \subset \partial K}} h_F^3
  \| \nabla \widetilde v_{h,a}(t)\cdot \nu \|_{L^2(F)}^2.
\end{align*}
The total local contribution is defined as $\mathcal{E}^n_{K} = \mathcal{E}^n_{1, K} + \mathcal{E}^n_{2, K} + \mathcal{E}^n_{3, K}$. Accounting for the full spatial and temporal domain gives rise to the partial estimators 
\[
  \mathcal{E}_j \defeq \sum_{n =0}^{N-1} \sum_{K \in \mathcal{T}_h} \mathcal{E}^n_{j,K} \qquad \text{for }j \in\{1,2,3\},
\]
and using the notation 
\[
  \mathcal{E}_0 \defeq \|\rho(0) - \bar \rho(0)\|_{H^{-1}(\Omega)}^2
\]
a decomposition of the a posteriori bound is given by
\begin{equation}
  \mathcal{E} = \left(1+\frac{1}{\gamma_{0,T}}\right)e^{\Lambda_T T}\Biggl( \mathcal{E}_0 + \mathcal{E}_1 + \mathcal{E}_2 + \mathcal{E}_3 \Biggr).
\end{equation}
In the following experiments, we study the individual estimator components $\mathcal E_j$, $j\in\{1,2,3\}$, and their sum. We omit the initial error $\mathcal E_0$ and the factor $(1+\gamma_{0,T}^{-1})e^{\Lambda_T T}$ from these comparisons. For the cell integrals, we use two-point Gauss--Legendre quadrature in 1D. In 2D, we use the tensor-product auxiliary mesh $\mathcal T_h^\dag$ obtained by subdividing each Cartesian cell into four equal squares. On this mesh, we compute the flux potential reconstruction of Remark~\ref{rem:tensor-product-auxiliary-meshes} and the auxiliary Galerkin approximation of Remark~\ref{rem:cartesian-auxiliary-elliptic}, and we use midpoint quadrature on each subsquare.

\subsection{Optimal scaling in 1D}\label{sec:optimality-1D}
In the first numerical experiment, we illustrate the behavior of the error estimator in 1D under global mesh refinement and demonstrate its optimality. We consider the computational domain $\Omega=(0,1)$ and choose the parameters
\[
  \alpha = 5000, \quad \beta= 800
  \]
  along with the sensitivity function and initial datum
  \[
    \gamma(v) = e^{-0.2v}, \quad \rho_0(x) = \sin^2(\pi x).
  \]

  \begin{figure}[t]
     \iftoggle{usetikz}{\input{figs/convergence-experiment-1D.tikz}}{\includegraphics{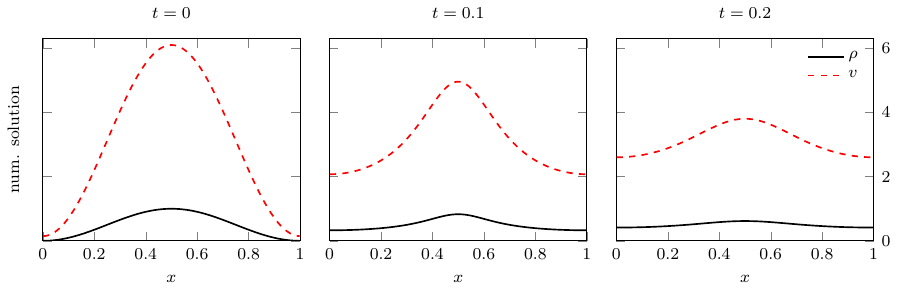}}
    \caption{Cell and attractant density over the domain $\Omega=(0,1)$ at three time instants in the first numerical experiment computed by the finite volume scheme.}\label{fig:experiment1}
  \end{figure}

  In the numerical solution shown in Fig.~\ref{fig:experiment1} both cells and chemoattractant slightly aggregate in the center of the domain after a short time ($t=0.1$), but overall diffusive effects dominate, leading to smooth spatial profiles at later times ($t=0.2$).

  \begin{table}[tb]
    \centering
    \caption{Components of the error estimator under global mesh refinement in the numerical experiment shown in Fig.~\ref{fig:experiment1}.\\[5pt]}\label{tab:convergence}
    \begin{tabular}{@{}rllllll@{}}
      \toprule
      $|\mathcal T_h|$ & $\mathcal{E}_1$ & EoC & $\mathcal{E}_2$ & EoC & $\mathcal{E}_3$ & EoC \\
      \midrule
      25 & $1.145 \times 10^{-4}$ &  & $3.202 \times 10^{-3}$ &  & $6.072$ &  \\
      50 & $2.892 \times 10^{-5}$ & 1.98 & $8.012 \times 10^{-4}$ & 2.00 & $1.281$ & 2.25 \\
      100 & $7.250 \times 10^{-6}$ & 2.00 & $2.004 \times 10^{-4}$ & 2.00 & $3.135 \times 10^{-1}$ & 2.03 \\
      200 & $1.814 \times 10^{-6}$ & 2.00 & $5.012 \times 10^{-5}$ & 2.00 & $7.807 \times 10^{-2}$ & 2.01 \\
      400 & $4.535 \times 10^{-7}$ & 2.00 & $1.255 \times 10^{-5}$ & 2.00 & $1.950 \times 10^{-2}$ & 2.00 \\
      800 & $1.134 \times 10^{-7}$ & 2.00 & $3.159 \times 10^{-6}$ & 1.99 & $4.874 \times 10^{-3}$ & 2.00 \\
      \bottomrule
    \end{tabular}
  \end{table}
  
  In the scheme we fix the small time increment $\Delta t = 10^{-5}$ to ensure that spatial errors dominate the temporal ones and compute numerical solutions at the final time $T=0.2$ on uniform grids with successively increasing numbers of grid cells. The components of the estimator for successively refined grid resolutions are presented in Table~\ref{tab:convergence} along with their estimated order of convergence (EoC)\footnote{The EoC is computed by the formula $EoC= \log2(\mathcal{E}_j^1/\mathcal{E}_j^2)$ with $\mathcal{E}_j^1$ and $\mathcal{E}_j^2$ denoting the partial estimators in two consecutive rows of the table.}. We observe that the partial estimator $\mathcal{E}_3$ accounting for the discretization error of the attractant dominates the other parts of the estimator in magnitude. All partial estimators decrease as the grid is refined, and our computations clearly point to second order convergence for all three components indicating optimality of the estimator.

\subsection{Local estimators and adaptive mesh refinement in 1D}\label{sec:refinementexperiment}
In this section, we investigate the localized estimators in another 1D experiment adopted from~\cite{2025HerdaTrescasesZurek} and embed them in a new adaptive mesh refinement method. The computational domain is again chosen as $\Omega=(0,1)$ with the parameters $\alpha=1000$ and $\beta=100$, and the sensitivity function and initial datum are given by
\[
  \gamma(v) = e^{-v}\quad \text{and} \quad \rho_0(v) = 15 x^2 (1- x)^2.
\]
\begin{figure}
   \iftoggle{usetikz}{\input{figs/localized-estimators-refinement.tikz}}{\includegraphics{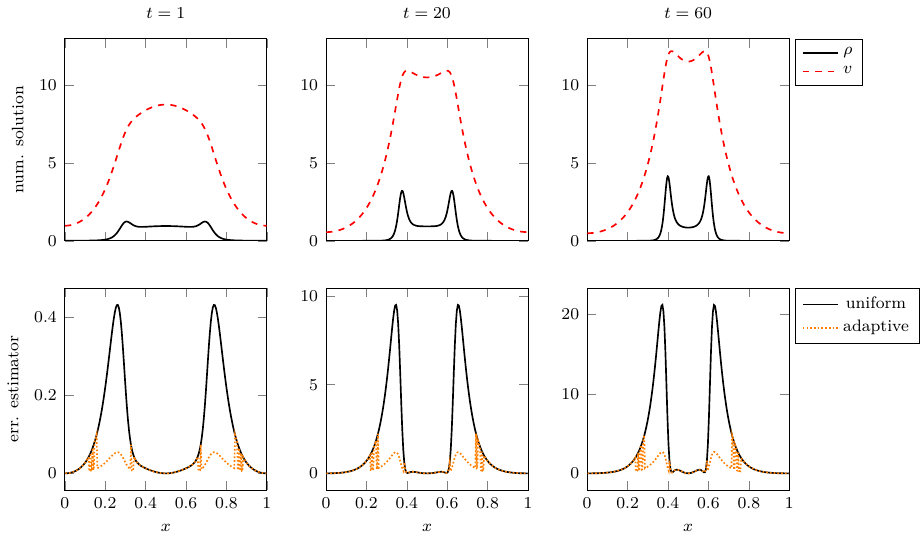}}
  \caption{Cell and attractant density over the domain $\Omega=(0,1)$ at three time instants in the second numerical experiment (upper row) along with the local distribution of the error estimator in the employed uniform numerical scheme and an adaptive mesh refinement method (lower row). The cells segregate into two moving concentrations with the estimator capturing their position and steep gradients.}\label{fig:localized-estimators-refinement}
\end{figure}
We fix the time increment $\Delta t = 0.01$ and first compute a numerical solution over a uniform grid with 200 cells. The numerical solution is shown in Fig.~\ref{fig:localized-estimators-refinement} (upper row): the cells arrange in two segregated concentrations at the beginning of the simulation (around $t=1$), which then further aggregate and move towards each other with decreasing velocity (see time instants $t=20$ and $t=60$). The chemical attractant accumulates over time with its density forming a local minimum at the center of the domain.

We also present the spatial distribution of the total local contributions to the estimator from the numerical simulation on the uniform mesh in Fig.~\ref{fig:localized-estimators-refinement} (lower row). The largest error indicators appear at the position of the moving concentrations and around their outer boundaries; the maximum increases over time together with the observed increase of the concentrations.

To improve the numerical accuracy, adaptive mesh refinement methods refine cells where the error estimator $\mathcal E_K^n$ is relatively large and coarsen cells where $\mathcal E_K^n$ is relatively small \cite{ohlberger2009}. As a proof of concept, we implement a mesh adaptation method based on the new error estimator following the approach in \cite{kolbeStudyTimeDiscretization2016, kolbe2022adaptrectanmesh}.
In short, we alternate between the solution update to the next time instant according to our scheme from Section 4 and the mesh update. During the latter, grid cells are refined via bisection to the next \emph{refinement level} if the \emph{monitor function} $\theta_K^n$ exceeds the threshold $0.5$ and coarsened via a merge to the previous refinement level if $\theta_K^n<0.25$. In the process, we use the monitor function based on the error estimator
\[
  \theta_K^n = \frac{\mathcal E_K^n}{\theta^n}, \qquad \theta^n = \max_{K \in \mathcal T_h^n} \mathcal E_K^n
  \]
and ensure that the refinement level does not exceed the fixed maximum level $m$ and that the refinement levels of adjacent cells differ by at most one.

Starting from a uniform grid with 200 cells at refinement level 0, we employ the adaptive method to compute a numerical solution in the time interval $[0,60]$ initially allowing only for the maximal refinement level $m=1$. The adaptive method employs a time weighted average\footnote{Time weighted average numbers of grid cells are computed by the formula $\frac{1}{T} \sum_{n=0}^{N-1} \Delta t^n |\mathcal{T}_h^n|$.} of $257.87$ grid cells and the spatial distribution of the error estimator for this method is indicated in Fig.~\ref{fig:localized-estimators-refinement} (lower row). The adaptive method significantly reduces the estimator in critical regions, where large values were observed in the uniform method.

\begin{figure}
  \iftoggle{usetikz}{\input{figs/refinement-estimator-over-time.tikz}}{\includegraphics{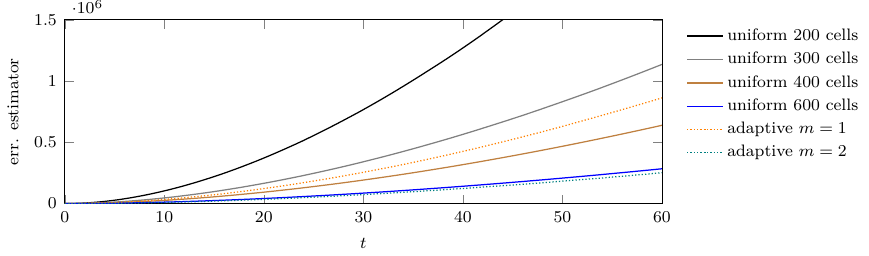}}
  \caption{The error estimator $\mathcal E_1 + \mathcal E_2 + \mathcal E_3$ over the simulation time interval $[0, 60]$ for different uniform and adaptive methods. The adaptive methods with maximum refinement level $m=1$ (time weighted average of $257.87$ grid cells) and $m=2$ (time weighted average of $416.97$ grid cells) outperform the uniform methods employing comparable numbers of grid cells.}\label{fig:refinement-estimator-over-time}
\end{figure}

For a detailed comparison, we present the error estimator $\mathcal E_1 + \mathcal E_2 + \mathcal E_3$ over the relevant time interval for different uniform and adaptive simulations in Fig.~\ref{fig:refinement-estimator-over-time}. Although the discussed adaptive method with maximal refinement level $m=1$ employs fewer cells on average, the corresponding error estimator stays below the error estimator of the uniform method on 300 grid cells over the full time interval. The advantage of the adaptive strategy becomes even more pronounced when allowing for the maximal refinement level $m=2$: even though the time weighted average is only $416.97$ grid cells, the adaptive method outperforms the uniform method on 600 cells in terms of the error estimator.

To evaluate the efficiency of the adaptive method we computed the runtime of the methods compared in Fig.~\ref{fig:refinement-estimator-over-time} up to the final time $T=60$ taking the median of 10 sample computations in each case\footnote{Computations were run on an M4 MacBook Pro with 24 GB of memory running macOS 27. Relative times were computed with respect to the uniform method on 200 cells.}. We found that the uniform methods, requiring neither an evaluation of the estimator nor a mesh update, had relative runtimes of $1$, $2.11$, $3.59$ and $10.42$ over $200$, $300$, $400$ and $600$ cells. The adaptive method had a relative runtime of $4.04$ for $m=1$ and $6.44$ for $m=2$. Given the error relations in Fig.~\ref{fig:refinement-estimator-over-time} the runtime computations show that adaptation limited by the low maximal level $m=1$ is outperformed in runtime by uniform methods yielding comparable errors. However, allowing for the maximum level $m=2$ we observe a significant gain in efficiency of the adaptive method (comparing the uniform method on 600 cells and the adaptive method with $m=2$).        

\subsection{Optimal scaling in 2D}
In this section, we validate the optimal scaling behavior of our estimator $\mathcal E$ in a numerical experiment in two space dimensions. To this end, we fix the square domain $\Omega=(0,1)^2$ and the production and decay parameters $\alpha=21000$ and $\beta=500$. As in Section~\ref{sec:refinementexperiment}, we assume that the sensitivity function is given by $\gamma(v)=e^{-v}$, and we select the initial cell density
\[
  \rho_0(\mathbf x) = \frac{\pi^3}{400} \sin^2(\pi x) \sin^2(\pi y).
  \]
\begin{figure}
   \iftoggle{usetikz}{\input{figs/2D-experiment.tikz}}{\includegraphics[width=\textwidth]{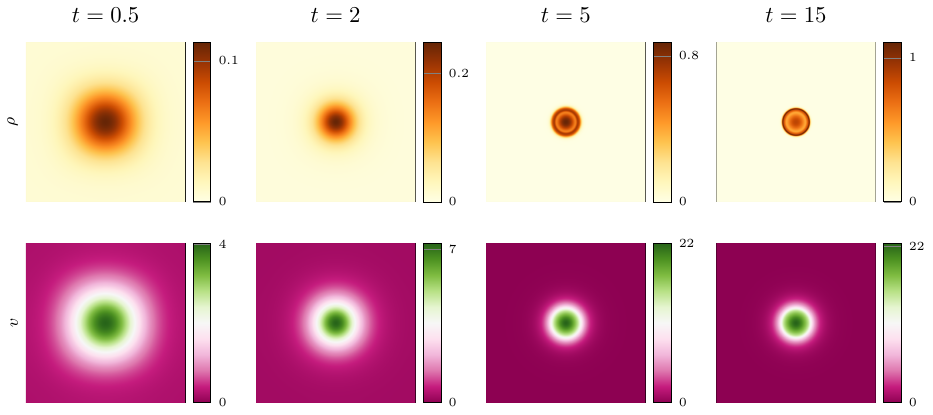}}
  \caption{Cell density (upper row) and attractant concentration (lower row) at four time instants in the 2D experiment over the domain $\Omega=(0,1)^2$. The numerical solution was computed by the finite volume scheme on a mesh using $400 \times 400$ cells. The cells aggregate and form ring-like patterns.}\label{fig:2D-experiment}
\end{figure}
Fig.~\ref{fig:2D-experiment} showcases the spatio-temporal dynamics of the cell density and the attractant concentration: the cells aggregate in the center of the domain (around $t=2$), where they arrange in ring-like patterns ($t=5$ and $t=15$). The chemoattractant concentrates at the cells' positions and seems to reach a steady state (at $t=5$) with high concentration ($v\approx 22$) at the center.

  As in Section~\ref{sec:optimality-1D}, we study the behavior of the partial estimators under global mesh refinement. Therefore, in this experiment we fix the final time $T=0.5$ and the small time increment $\Delta t= 10^{-4}$. Numerical solutions are computed on successively refined uniform grids, and the corresponding partial estimators are presented in Table~\ref{tab:convergence-2D} together with their EoC.
  Also in the case of two space dimensions, the partial estimators decay as the mesh is refined, and as the EoC approaches 2 in all three cases, the numerical results confirm the optimality of the estimator~\eqref{eq:computable-estimate}.
  
 \begin{table}
    \centering
    \caption{Components of the error estimator under global mesh refinement in the 2D experiment shown in Fig.~\ref{fig:2D-experiment}.\\[5pt]}\label{tab:convergence-2D}
    \begin{tabular}{@{}rllllll@{}}
      \toprule
      $\mathcal T_h$ & $\mathcal{E}_1$ & EoC & $\mathcal{E}_2$ & EoC & $\mathcal{E}_3$ & EoC \\
      \midrule
      $5 \times 5$ & $1.006 \times 10^{-5}$ &  & $4.098 \times 10^{-5}$ &  & $417.9$ &  \\
      $10 \times 10$ & $5.979 \times 10^{-6}$ & 0.75 & $9.145 \times 10^{-6}$ & 2.16 & $45.74$ & 3.19 \\
      $20 \times 20$ & $1.517 \times 10^{-6}$ & 1.98 & $1.796 \times 10^{-6}$ & 2.35 & $8.274$ & 2.47 \\
      $40 \times 40$ & $4.906 \times 10^{-7}$ & 1.63 & $4.495 \times 10^{-7}$ & 2.00 & $1.110$ & 2.90 \\
      $80 \times 80$ & $1.267 \times 10^{-7}$ & 1.95 & $1.123 \times 10^{-7}$ & 2.00 & $1.718 \times 10^{-1}$ & 2.69 \\
      $160 \times 160$ & $3.199 \times 10^{-8}$ & 1.99 & $2.857 \times 10^{-8}$ & 1.97 & $3.661 \times 10^{-2}$ & 2.23 \\
      \bottomrule
    \end{tabular}
  \end{table}

  \section{Conclusion}
  We have derived computable a posteriori error estimates for a finite volume approximation of the parabolic--elliptic local sensing Keller--Segel model. By establishing a stability estimate that requires no smallness condition on the residuals, we obtained upper bounds for the cell density error in the $L^\infty(0, T; H^{-1}(\Omega))$ and $L^2(0, T; L^2(\Omega))$ norms. The resulting estimator, which utilizes reconstructions of the flux, flux potential, and chemoattractant concentration, demonstrated optimal scaling in numerical experiments in one and two space dimensions. Furthermore, we proposed an adaptive mesh refinement method steered by the estimator and verified its efficiency in the one-dimensional case, showing that it can significantly reduce the error while employing fewer grid cells than uniform meshes. Our results provide a framework for quantifying errors and optimizing computational resources when simulating complex aggregation patterns in local sensing chemotaxis models.

\paragraph{Acknowledgements}
  The authors thank Jan Giesselmann for his valuable comments and suggestions.

  \paragraph{Funding}
  NK thanks the Deutsche Forschungsgemeinschaft (DFG, German Research Foundation) for the financial support through SPP 2311 (Robust Coupling of Continuum-Bio\-me\-chanical In Silico Models to Establish Active Biological System Models for Later Use in Clinical Applications – Co-Design of Modelling, Numerics and Usability) within Project~548864771 and through SPP 2410 (Hyperbolic Balance Laws in Fluid Mechanics: Complexity, Scales, Randomness) within Project~525842915. KK was supported by the National Research Foundation of Korea (NRF) grants funded by the Korea government (MSIT) (Nos. 00560003 and 2023R1A2C1003120).
  \bibliographystyle{abbrvurl} 
  \bibliography{references.bib}

\end{document}